\documentclass[hidelinks,onefignum,onetabnum]{siamart220329}

\usepackage{lipsum}
\usepackage{amsfonts}
\usepackage{graphicx}
\usepackage{epstopdf}
\usepackage{algorithmic}
\usepackage{amsmath}

\usepackage{amssymb}

\usepackage{bm}
\usepackage{dsfont}
\usepackage{subfig} % Numbered and caption subfigures using \subfloat
\usepackage{cite}
\usepackage{diagbox}
\usepackage{hyperref}

\ifpdf
  \DeclareGraphicsExtensions{.eps,.pdf,.png,.jpg}
\else
  \DeclareGraphicsExtensions{.eps}
\fi
\newcommand{\mathbbm}[1]{\text{\usefont{U}{bbm}{m}{n}#1}} % From mathbbm.sty
\newcommand{\ICNN}{\text{ICNN}_{\bm{\theta}}}

\newsiamremark{remark}{Remark}
\newsiamremark{hypothesis}{Hypothesis}
\crefname{hypothesis}{Hypothesis}{Hypotheses}
\newsiamthm{claim}{Claim}

\headers{Structure-preserving NN in data-driven rheology}{N. Parolini, A. Poiatti, M. Verani, and J. Ven\'e}

\title{Structure-preserving neural networks in data-driven rheological models
}

\author{
Nicola Parolini\thanks{MOX, Dipartimento di Matematica, Politecnico di Milano, Piazza Leonardo da Vinci 32 - 20133 Milano, Italy
  (\email{nicola.parolini@polimi.it}, \email{julian.vene@polimi.it}, \email{marco.verani@polimi.it}).}
\and 
Andrea Poiatti\thanks{Dipartimento di Matematica, Politecnico di Milano, Piazza Leonardo da Vinci 32 - 20133 Milano, Italy
  (\email{andrea.poiatti@polimi.it}).}
\and 
Julian Vené\footnotemark[2]
\and 
Marco Verani\footnotemark[2]
}

\usepackage{amsopn}

\begin{document}

\maketitle

% REQUIRED
\begin{abstract}
In this paper we address the importance and the impact of employing structure preserving neural networks as surrogate of the analytical  physics-based models typically employed to describe the rheology of non-Newtonian fluids in Stokes flows. In particular, we propose and test on real-world scenarios a novel strategy to build data-driven rheological models based on the use of Input-Output Convex Neural Networks (ICNNs), a special class of feedforward neural network scalar valued functions that are convex with respect to their inputs. Moreover, we show, through a detailed campaign of numerical experiments,  that the use of ICNNs is of paramount importance to guarantee the well-posedness of the associated non-Newtonian Stokes differential problem. Finally, building upon a novel perturbation result for non-Newtonian Stokes  problems, we study the impact of our data-driven ICNN based rheological model on the accuracy of the finite element approximation.
\end{abstract}

% REQUIRED

% REQUIRED

\section{Introduction}
\label{sec:introduction}
Rheological models are closed-form mathematical expressions describing the relation between the stress and the rate of deformation of the fluid. The study of the rheology of complex fluids has been a subject of interest for more than one century: based on the experimental evidence showing different types of rheological responses, many constitutive models have been proposed.
Generally speaking, there are two main approaches for obtaining shear–stress relationships: data-driven modeling based on empirical observations; and physics-based models derived from first principles based on the material's underlying structure.
In the context of the second approach, even limiting the discussion to the so-called \textit{generalized Newtonian}, in which the stress is assumed to be only a function of shear rate and does not depend upon the history of deformation, several  rheological models have been put forward and are routinely adopted in engineering applications, ranging from the simple power-law model, in which the shear stress is proportional to a power of the shear rate, to more involved models such as Cross \cite{cross1965}, Carreau and Carreau-Yasuda \cite{carreau} models, in which the power law dependence is limited to intermediate shear rates, while asymptotic constant viscosities are achieved for both low and high shear rates.                         Clearly, the use of physics-based rheological models requires an {\em a-priori} choice of the analytical dependence between shear rate and viscosity, leaving to the calibration phase the choice of the parameters to best fit the experimental data. 

To overcome the need of the a-priori (to some extent arbitrary) choice of the model, which may be limiting when complex flow behaviors are considered, data-driven rheological models can represent a suitable and valid alternative. In recent years, constitutive models based on neural networks (NN) have received growing attention (see, e.g., \cite{Reyes2021,Chen2022,As'ad20222738,Klein2022,HUANG2022104856,Howard2023507,Karniadakis:2023} and the references therein), owing to their excellent nonlinear function fitting capabilities, while remaining form-agnostic, which makes them desirable for a general shear–stress modeling framework. Moreover, NNs possess two computational attractive features: (1) a rapid online execution time, which is particularly beneficial when the evaluation of a constitutive law can dominate the simulation cost; (2) the current availability and efficiency of NN automatic differentiation packages, which simplifies the task of obtaining tangent problems to solve the nonlinear system stemming, e.g., from the Finite Element (FE) discretization of the associated differential problem. For the solution of algebraic problems employing machine learning techniques see, e.g., \cite{klawonn:2021} and the references therein.

However, when data-driven models are employed in combination with systems of partial differential equations (PDEs), there are key issues that have to be properly addressed to make this approach  amenable to scientific simulations. Indeed, the data-driven model not only must satisfy physical assumptions, but also mathematical properties to ensure the existence (and potentially, the uniqueness) of  solutions of the associated differential problem.
Thus, it becomes of paramount importance  the integration of {\em structure preserving} neural networks into  partial differential equations (PDEs) that guarantee the well-posedness of the mathematical problem.

 \subsection{Our contribution}
 In this paper we consider the following prototypical non-Newtonian Stokes problem:
  \begin{equation} \label{eq:intro}
    \begin{cases}
        - \nabla \cdot [\tau(\varepsilon(\mathbf{u}))] + \nabla p = 
        \mathbf{f} \qquad \text{in } \Omega,  \\
        \nabla \cdot \mathbf{u} = 0 \qquad \text{in } \Omega, \\
        u =0 \qquad \text{on } \Omega, \\
    \end{cases}
\end{equation}
where we assume that the stress tensor $\tau$ has to be deduced from experimental data. 
It is worth noticing that from a mathematical point of view (see, e.g., \cite{barretnonnewton}) the stress tensor $\tau$ must satisfy a certain set of conditions (cf. Assumptions {\bf (A)} below) to ensure the well-posedness of the differential problem \eqref{eq:intro}.   In view of this remark,  a  common practice to find $\tau$ and the associated solution $(\mathbf{u},p)$ is the following: 
\begin{itemize}
\item {\it a priori} select  a law for $\tau$ that is known to ensure the well-posedness of the problem (e.g.  Carreau or power law, that typically depend on a certain set of parameters);
\item employ experimental data to fit the parameters of the selected law;
\item numerically solve \eqref{eq:intro}.
 \end{itemize}
 As mentioned earlier, the potential difficulty of this approach is the {\it a priori} choice of the law for $\tau$,  which   can be (to some extent) arbitrary, since it is not adapted to data.
 For this reason,  in this paper we consider a data-driven discovery of the stress tensor $\tau$ from experimental data, based on neural networks.  In particular,  to ensure the well-posedness of the differential problem we make the following choice 
 \begin{equation*}
   \tau_{\bm{\theta}}(\varepsilon(\mathbf{u})) = \pm\text{ICNN}_{\bm{\theta}}( |\varepsilon(\mathbf{u}) |)\varepsilon(\mathbf{u})
\end{equation*}
 where $\ICNN$: $\mathbb{R}\to\mathbb{R}$ is an Input-Convex neural network (see \cite{amos2017input,sivaprasad2021curious}) whose parameters ${\bm{\theta}}$ are trained on the specific set of experimental data and the sign is automatically learned depending on the nature of fluid (shear thinning or shear thickening). 
 Differently from standard feed-forward neural networks,  Input-Convex Neural Network functions  are convex functions with respect to the input variable. In view of this crucial property,  it turns out that  the following problem (and its discrete approximation) is well posed:
 \begin{equation} \label{eq:intro:2}
    \begin{cases}
        - \nabla \cdot [\tau_{{\bm{\theta}}}(\varepsilon(\mathbf{u}))] + \nabla p = 
        \mathbf{f} \qquad \text{in } \Omega,  \\
        \nabla \cdot \mathbf{u} = 0 \qquad \text{in } \Omega, \\
        u =0 \qquad \text{on } \Omega.\\	
    \end{cases}
\end{equation}    
To conclude, we mention that in \cite{Reyes2021}
a data-driven neural network based approach to learn viscosity models of two non-Newtonian systems (polymer melts and suspensions of particles) using only velocity measurements has been studied by employing Physics Informed Neural Networks \cite{Raissi2019686}. However, in \cite{Reyes2021} standard neural network have been employed and the crucial issue of the well-posedness of the resulting differential problem \eqref{eq:intro:2} is not addressed.
 
Throughout the paper we will use the notation $x\lesssim y$ with the meaning $x\leq c y$, with $c$ positive constant independent of the discretization parameters.      
\subsection{Outline}
In Section \ref{sec:FEapprox}, we introduce the  Stokes equation for non-Newtonian fluids   and recall useful theoretical results concerning the well-posedness of the continuous problem and the approximation properties of its finite element approximation; these latter properties are exemplified through a set of numerical experiments, that are instrumental for the subsequent discussion.
In Section \ref{sec:ICNN}, we briefly recall Input-Convex Neural Network (ICNN) scalar valued functions and we show how to employ them to build  data-driven rheological models that are conformal with the well-posedness of the differential problem.
In Section \ref{sec:ICNNStokes}, we address the finite element approximation of the non-Newtonian Stokes problem governed by the data-driven ICNN rheological model. In particular, building upon a perturbation-type result for non-Newtonian Stokes problem, we connect the approximation properties of the ICNN and the ones of the finite element method. Finally, in Section  
\ref{S:Concl} we draw our conclusions.

%%%%%%%%%%%%%%%%%%%%5
\section{Non-Newtonian Stokes flows}
\label{sec:FEapprox}
In this section we recall the continuous problem together with its weak formulation and highlight suitable hypotheses guaranteeing existence and uniqueness of the solution (Section \ref{sec:FEapprox:1}). Later, we introduce the finite element discretization (Section \ref{sec:FEapprox:2}), whose approximation properties are exemplified and discussed through a set of numerical tests (Section \ref{sec:numerics_1}).  
\subsection{Continuous problem and weak formulation}
\label{sec:FEapprox:1}
Let us consider an incompressible fluid inside a domain $\Omega \subset \mathbbm{R}^2$, with Lipschitz boundary $\partial \Omega$. We denote
by $\mathbf{u}: \Omega \rightarrow \mathbbm{R}^2$ the velocity field and by $p: \Omega \rightarrow \mathbbm{R}$ the pressure. We suppose that $(\mathbf{u},p)$ satisfies the following equation:
\begin{equation} \label{eq:s_ns}
    \begin{cases}
        - \nabla \cdot [\tau(x,\varepsilon(\mathbf{u}))] + \nabla p = 
        \mathbf{f} &\qquad \text{in } \Omega,  \\
        \nabla \cdot \mathbf{u} = 0 &\qquad \text{in } \Omega, \\
    \end{cases}
\end{equation}
coupled with the homogeneous Dirichlet boundary condition for $\mathbf{u}$ and zero mean condition for the pressure $p$:
\begin{equation*}
    \mathbf{u} = \mathbf{0} \text{ on } \partial \Omega, \qquad \int_\Omega p \, dx = 0.
\end{equation*}
Here $\tau$ denotes the stress tensor, which is a suitable function of the symmetric strain rate tensor $\varepsilon(\mathbf{u})$ defined as:
\begin{equation*}
    \varepsilon(\mathbf{u}) := \dfrac{1}{2}(\nabla \mathbf{u} + \nabla \mathbf{u}^T)
\end{equation*}
and the term $\mathbf{f}$ is a given body force.
In the framework, adopted in this paper, of the so-called {\em generalized Newtonian fluids}, the stress tensor is assumed to be a function of the shear rate only, i.e. $\tau(x,\varepsilon(\mathbf{u})) = \tau(\varepsilon(\mathbf{u}))$ and does not depend on the history of deformation. More precisely, we have:
\begin{equation*}
    \tau(\varepsilon(\mathbf{u})) = k(\left\lvert \varepsilon(\mathbf{u})\right\rvert)\varepsilon(\mathbf{u})
\end{equation*}
where the function $k:\mathbb{R}\to \mathbb{R}$ represents the viscosity of the fluid, which results to be a constant in case of a Newtonian fluid, while the norm $\left\lvert \cdot \right\vert^2 $ represents the squared Frobenius norm, i.e. for $\mathbf{K} \in \mathbbm{R}^{n \times n}$ real matrix, $\left\lvert \mathbf{K} \right\rvert^2 =  \sum_{i,j=1}^n K_{ij}^2$. Popular choices for $k(\cdot)$ are the Carreau law:
    \begin{equation}\label{eqn:carreau}
        k(t) = k_\infty + (k_0 - k_\infty)(1+ \lambda t^2)^{(n-2)/2},
    \end{equation} 
    and the power law:
    
    \begin{equation}\label{eqn:power}
        k(t) = k_0 t^{n-2}
    \end{equation}
where $k_0 > k_\infty \geq 0$, $\lambda > 0$ and $n \in (1,2)$ for pseudo-plastic fluid and $n>2$ for a dilatant fluid. The case $n=2$ corresponds to Newtonian fluids.

About the mathematical analysis of problem \eqref{eq:s_ns}, the literature is quite flourishing: we only refer to the works \cite{N1,N2,Malek,B1,B2,barretnonnewton} and references therein (see also \cite{Bernardi,GPPV} and their references for more recent contributions on the topic, under more general assumptions of non-isothermal fluids). In particular, in this paper we will adopt the assumptions (and the results) presented in \cite{barretnonnewton}. We thus  introduce a set of  assumptions for $k(\cdot)$ which guarantee the existence and uniqueness of the (weak) solution to \eqref{eq:s_ns}: \\

\textbf{Assumptions (A)}: We assume that $k \in C(0,\infty)$ and that there exist constants  $r \in (1, \infty)$, $\alpha \in [0,1]$ and 
$\varepsilon, C, M > 0$ such that:
\begin{equation*}\label{ass_1}
%    \begin{multlined}
    k(t) \leq C [ t^\alpha(1+t)^{1-\alpha}    ]^{r-2} \qquad \forall t \geq 0,
%    \end{multlined}
\end{equation*}
\begin{equation*}\label{ass_2}
%\begin{multlined}
    |k(t)t- k(s)s| \leq  C|t-s|[(t+s)^\alpha(1+t+s)^{1-\alpha}]^{r-2}  \qquad \forall t,s > 0 \text{ such that } |s/t -1| \leq \varepsilon,
%\end{multlined}
\end{equation*}

\begin{equation*}\label{ass_3}
%\begin{multlined}
    k(t)t- k(s)s \geq M(t-s)[(t+s)^\alpha(1+t+s)^{1-\alpha}]^{r-2}\qquad  \forall t \geq s \geq 0.
%\end{multlined}
\end{equation*}
\begin{remark}\label{Rem:assumptions}
The Assumptions \textbf{(A)} are satisfied by the Carreau law \eqref{eqn:carreau} with $\alpha=0$ and $r=n$ if $k_\infty=0$ and $r=2$ if $ n\in (1,2]$ and $k_\infty>0 $ and by the power law \eqref{eqn:power} 
with $\alpha=1$ and $r=n$
(cf. \cite[Remark 2.1]{barretnonnewton}).
The parameter $\alpha$ entering in the Assumptions \textbf{(A)} measures the degree of degeneracy of $k(\cdot)$ for a given $r$: the closer $\alpha$ is to one the more degenerate $k(\cdot)$ is.
\end{remark}

 We are now ready to introduce the weak formulation of (\ref{eq:s_ns}). Following \cite{barretnonnewton} let us assume that $k(t)$ satisfies  Assumptions \textbf{(A)} and for $\mathbf{v} \in [W_0^{1,r}(\Omega)]^2$ we define:
\begin{equation}
    J(\mathbf{v}) := \int_\Omega \Big[   \int_0^{|\varepsilon(\mathbf{v}) |} k(t)tdt \Big] - \left\langle \mathbf{f},\mathbf{v}\right\rangle,
\end{equation}
where $\langle f,g\rangle:=\int_\Omega fg\ dx$. 
We set $\mathbf{X} = [W_0^{1,r}(\Omega)]^2$, and let $\left\langle \cdot, \cdot \right\rangle$ be the duality pairing between the dual space $\mathbf{X}^*$ and $\mathbf{X}$. It is easy to check that 
$J(\cdot)$ is Gateaux differentiable on $\mathbf{X}$ with:
\begin{equation*}
%\begin{multlined}
    \left\langle J'(\mathbf{w}_1), \mathbf{w}_2 \right\rangle_{\mathbf{X}^*} = \left\langle A\mathbf{w}_1, \mathbf{w}_2 \right\rangle_{\mathbf{X}^*} - 
    \left\langle \mathbf{f}, \mathbf{w}_2\right\rangle \qquad \forall \mathbf{w}_1,\mathbf{w}_2 \in \mathbf{X},
%\end{multlined}
\end{equation*}
where $A:\mathbf{X} \rightarrow \mathbf{X}^*$ is such that
\begin{equation*}
    \left\langle A\mathbf{w}_1, \mathbf{w}_2 \right\rangle_{\mathbf{X}^*} =  \left\langle k(\left\vert \varepsilon(\mathbf{w}_1)\right\vert)\varepsilon(\mathbf{w}_1),\varepsilon(\mathbf{w}_2)\right\rangle. 
\end{equation*}

In \cite{barretnonnewton} it is proved that  if $k(\cdot)$ satisfies the assumptions \textbf{(A)} then the functional $J(\mathbf{v})$ is strictly convex on  $\mathbf{X}$ and the following result holds:
    \begin{proposition}
    Let $\mathbf{V} := \{\mathbf{v} \in \mathbf{X} : \nabla \cdot \mathbf{v} = 0 \text{ in } \Omega \}$. The problem:
    \begin{equation*}
        J(\mathbf{u}) \leq J(\mathbf{v}) \qquad \forall \mathbf{v} \in \mathbf{V}
    \end{equation*}
    or equivalently:
    \begin{equation}\label{eqn:R_continous}
        \left\langle A\mathbf{u},\mathbf{v}\right\rangle_{\mathbf{X}^*} = \left\langle \mathbf{f},\mathbf{v}\right\rangle \qquad \forall \mathbf{v} \in \mathbf{V}
    \end{equation}
    admits a unique solution.
    \end{proposition}\label{prop:well-posed}

Let us define $M = L_0^{r'}(\Omega)$, where $r'$ is the conjugate exponent of $r$.  Considering  the mixed formulation of \eqref{eqn:R_continous} provides the weak formulation of (\ref{eq:s_ns}) which reads reads as follows: 
find $(\mathbf{u},p) \in \mathbf{X} \times M$ such that:
\begin{equation} \label{eqn:weak_stokes}
    \begin{cases}
        \left\langle A\mathbf{u}, \mathbf{w}\right\rangle_{\mathbf{X}^*} - \left\langle p, \nabla \cdot \mathbf{w}\right\rangle  = 
        \left\langle \mathbf{f},\mathbf{w}\right\rangle  \qquad \forall \mathbf{w} \in \mathbf{X},\\
        \left\langle \nabla \cdot \mathbf{u}, q   \right\rangle  = 0 \qquad \forall q \in M.\\
    \end{cases}
\end{equation}
For the well-posedness of \eqref{eqn:weak_stokes} we require the inf-sup condition:
\begin{equation}\label{eqn:babuska-brezzi-cond}
    \inf_{q \in M} \sup_{\mathbf{w} \in \mathbf{X}}\dfrac{\left\langle q, \nabla \cdot \mathbf{w} \right\rangle }{\|q\|_M \|\mathbf{w}\|_\mathbf{X}} \geq \beta > 0
\end{equation}
to hold.  In particular,  Amrouche and Girault (1990) proved in \cite{Amrouche1994} that when taking $\mathbf{X} = [W_0^{1,r}(\Omega)]^2$ and $ M = L_0^{r'}(\Omega)$, with $r$ and $r'$ conjugate exponents, there exist a constant $\beta(r)$ such that:
\begin{equation}\label{inf-sup:cont}
    \inf_{q \in M} \sup_{\mathbf{w} \in \mathbf{X}}\dfrac{\left\langle q, \nabla \cdot \mathbf{w} \right\rangle }{\|q\|_M \|\mathbf{w}\|_\mathbf{X}} \geq \beta(r) > 0.
\end{equation}
Hence (\ref{eqn:babuska-brezzi-cond}) holds, and therefore the existence of 
a unique solution $(\mathbf{u},p)$ in $(\ref{eqn:weak_stokes})$ is implied.  Note that  the unique solution $\mathbf{u}$ coincides with the solution of \eqref{eqn:R_continous}, as it can be seen from \eqref{eqn:weak_stokes} by restricting the test functions $\mathbf{w}$ to $\mathbf{V} \subset \mathbf{X}$.

\subsection{Finite Element approximation}
\label{sec:FEapprox:2}
In this section,  we briefly introduce the finite element approximation of \eqref{eqn:weak_stokes}.  Let $\mathbf{X}_h$ and $M_h$ be two finite dimensional spaces such that
\begin{equation*}
    \mathbf{X}_h \subset \mathbf{X} \cap [W^{1,\infty}(\Omega)]^2 \quad , \quad M_h \subset M \cap L^\infty(\Omega).
\end{equation*}
The discretized version of \eqref{eqn:weak_stokes} reads: 
find $(\mathbf{u}_h,p_h) \in \mathbf{X}_h \times M_h$ such that
\begin{equation} \label{eqn:mixed_weak_stokes_disc}
    \begin{cases}
        \left\langle A \mathbf{u}_h, \mathbf{w}_h \right\rangle_{\mathbf{X}^*} - \left\langle p_h, \nabla \cdot \mathbf{w}_h \right\rangle  = \left\langle \mathbf{f}, \mathbf{w}_h \right\rangle,
\\
        \left\langle \nabla \cdot \mathbf{u}_h, q_h   \right\rangle  = 0 \\
    \end{cases}
\end{equation}
holds $\forall \mathbf{w}_h \in \mathbf{X}_h, \forall q_h \in M_h$.

We further introduce the following three classical hypotheses (cf.  \cite{barretnonnewton}):
\begin{itemize}
    \item \textbf{(H1)} Approximation property of $\mathbf{X}_h$: there is a continuous linear operator $\pi_h:[W_0^{1,r}(\Omega)]^2 \rightarrow \mathbf{X}_h$ such that for 
        $ j = 0 ,...,m$ we have
        \begin{equation*}
%            \begin{multlined}
                \left\lVert \mathbf{w} - \pi_h \mathbf{w}\right\rVert_{[W^{j+1,r}(\Omega)]^2}  \leq 
        Ch^j \| \mathbf{w} \|_{[W^{j+1,r}(\Omega)]^2},
%            \end{multlined}
        \end{equation*}
        for any $\mathbf{w} \in [W_0^{1,r}(\Omega)]^2 \cap [W^{j+1,r}(\Omega)]^2$.
    \item \textbf{(H2)} Approximation property of $M_h$:
        there is a continuous linear operator $\rho_h: L^{r'}(\Omega) \rightarrow M_h$ such that for all $j= 0,...,m$ we have
        \begin{equation*}
%            \begin{multlined}
                \left\lVert q - \rho_h q\right\rVert_{L^{r'}(\Omega)}  \leq 
        Ch^j \| q \|_{L^{r'}(\Omega)}  \qquad \forall q \in L^{r'}(\Omega).
%            \end{multlined}
        \end{equation*}
          \item \textbf{(H3)} Discrete inf-sup condition: for any $r\in(1,\infty)$ there exists a constant $\beta_h(r)>0$ such that 
         \begin{equation*}
    \inf_{q_h \in M_h} \sup_{\mathbf{w_h} \in \mathbf{X_h}}\dfrac{\left\langle q_h, \nabla \cdot \mathbf{w_h} \right\rangle }{\|q_h\|_M \|\mathbf{w_h}\|_\mathbf{X}} \geq \beta_h(r) > 0.
\end{equation*}  
\end{itemize}
Under those assumptions,  the following result holds (cf.  \cite[Thm 4.1]{barretnonnewton})

\begin{theorem}\label{thm:convergence}
    Assume that $k(\cdot)$ satisfies Assumptions \textbf{(A)}. Let $(\mathbf{u},p)$ be the solution of ($\ref{eqn:weak_stokes}$) and let $(\mathbf{u}_h,p_h) \in \mathbf{X}_h \times M_h$ be the solution of 
\eqref{eqn:mixed_weak_stokes_disc}. Assume that \textbf{(H1)}-\textbf{(H2)} for $j = 1,...,m$ and \textbf{(H3)} hold. If $r \in (1,2]$, $\theta\in [r, 2+\alpha(r-2)]$ 
        \begin{equation}\label{eqn:err_uh}
%            \begin{multlined}
            \| \mathbf{u} - \mathbf{u}_h \|_{[W^{1,r}(\Omega)]^2} + \| p-p_h \|_{L^{\theta '}(\Omega)}^{\theta/[2(\theta-1)]} \leq C_1 h^{\theta j/2}
%               \end{multlined}
        \end{equation}
        where $
        C_1 = C( \left\lVert \mathbf{u} \right\rVert_{[W^{j+1,\theta}(\Omega)]^2}, \left\lVert p \right\rVert_{W^{j,r'}(\Omega)},$
        $ \beta_h(r)^{-1})$.

        If $r \in [2,\infty)$, $\theta\in [2+\alpha(r-2),r]$ we have 
        \begin{equation}\label{eqn:err_uh2}
%            \begin{multlined}
            \| \mathbf{u} - \mathbf{u}_h \|_{[W^{1,\theta}(\Omega)]^2} + \| p-p_h \|_{L^{r'(\Omega)}}^{2/\theta} \leq C_2 h^{j/(\theta-1)}
%               \end{multlined}
        \end{equation}
        where $
        C_2 = C( \left\lVert \mathbf{u} \right\rVert_{[W^{1,\infty}(\Omega)]^2},\left\lVert \mathbf{u} \right\rVert_{[W^{j+1,2}(\Omega)]^2}, \left\lVert p \right\rVert_{W^{j,\theta'}(\Omega)},$
        $ \beta_h(r)^{-1})$.
        
\end{theorem}
\begin{remark}\label{Rem:rate}
    If $k(\cdot)$ is not degenerate, i.e., $\alpha=0$ (cf. Remark \ref{Rem:assumptions}), then the error estimates in Theorem \ref{thm:convergence} are optimal on choosing $\theta=2$ under suitable regularity assumptions of the solution pair $({\bf u},p)$. The rates of convergence deteriorate as $k(\cdot)$ degenerates (i.e. $\alpha\to 1$). 
\end{remark}
\begin{remark}
As the study of the numerical approximation of \eqref{eq:s_ns} is not the focus of the paper, we limit our discussion to the pioneering Theorem \ref{thm:convergence} which is instrumental for the our purpose and nonetheless still represents nowadays a benchmark for the error analysis of the approximation of non-Newtonian Stokes problems.  We refer, e.g.,  to \cite{Hirn,R1,R2,R3,R4} and the references therein for more recent contributions in the field.
\end{remark}
\subsection{Numerical results}\label{sec:numerics_1}
In this section we show some numerical tests corroborating the error estimates  contained in Theorem \ref{thm:convergence}. The outcome of the simulations will be the base to study the impact of the use of neural networks as approximation of the stress tensor $\tau$ (see Section \ref{sec:ICNN} for the introduction of the particular class of employed neural networks and Section \ref{sec:ICNNStokes} for the impact on the numerical solution of \eqref{eq:s_ns}).  To this goal,  in the sequel  we consider the Carreau law \eqref{eqn:carreau} with parameters $k_0 = 2$, $k_\infty = 0$, $\lambda = 2$ and $n \in \{ 1.2, 1.6, 2, 2.4, 2.8\}$, inside the domain $\Omega = (-0.5,0.5)^2$. We recall (cf. Remark \ref{Rem:assumptions}) that in this case Assumptions ({\textbf{A}}) are satisfied with $r=n$ and $\alpha=0$. To evaluate the rate of the convergence of the approximation, the source term $\mathbf{f}$ is 
manufactured so that the exact solution $(\overline{\mathbf{u}},\overline{p})$ is explicitly given by:
\begin{align*}
    \overline{\mathbf{u}}(x,y) &=   \begin{bmatrix}
                            \phantom{-}5y \sin (x^2 + y^2) + 4y \sin(x^2 - y^2) \\
                            -5x \sin(x^2 + y^2) + 4x \sin(x^2 - y^2)
                        \end{bmatrix}\\
    \overline{p}(x,y) &= \sin(x+y).
\end{align*}
Dirichlet boundary conditions for velocity given by the exact solution are imposed on the domain boundary:
$
\mathbf{u}|_{\partial \Omega} = \overline{\mathbf{u}}|_{\partial \Omega}.
$ 
We consider the Taylor-Hood finite element spaces:
\begin{align*}
\mathbf{X}_h^j &:= \{\, \mathbf{w}_h \in [C(\Omega)]^2 : \mathbf{w}_h|_\tau \in P_j(\tau), \forall \tau \in T^h, 
              \text{and } \mathbf{w}_h = \overline{\mathbf{u}} \text{ on } \partial \Omega \}, \\
M_h^{j-1} &:= \{\, q_h \in L^{2}_0(\Omega) : q_h |_\tau \in P_{j-1}(\tau), \forall \tau \in T^h \},
\end{align*}
where $P_m$ is the space of all polynomials of degree less than or equal to $m$, and $T^h$ is a collection of disjoint open regular triangles $\tau$ such that $\bigcup_{\tau \in T^h}\tau = \Omega$. 

The numerical solution of the nonlinear system stemming from the discretization of \eqref{eqn:mixed_weak_stokes_disc} is carried out using the \textit{Newton method with trust region} \cite{NoceWrig06} available in Firedrake \cite{FiredrakeUserManual} through its PETSc backend \cite{balay1998petsc}.

In Table \ref{table:conv-rate}, we report the experimental convergence orders obtained for different values of the parameter $n$ and different polynomial degrees $j$. These results are in line with the theoretical predictions reported in Theorem \ref{thm:convergence} (see also Remark \ref{Rem:rate}). In particular, as expected when the Carreau model is considered, optimal convergence rates are achieved for both pseudo-plastic and dilatant fluids.

\begin{table}[H]
    \centering 
\begin{tabular}{|c|ccccc|ccccc|}
\hline
          & \multicolumn{5}{c|}{$W^{1,r}$ Velocity Error} & \multicolumn{5}{c|}{$L^{r^\prime}$ Pressure Error}  \\
\hline
\backslashbox{j}{r}         & 1.2&1.6 &2 &2.4 &2.8 &1.2&1.6 &2 &2.4 &2.8 \\
\hline
2 ($P_2/P_1$) & 2.00 & 2.00 & 2.00 & 2.00 & 2.00 & 2.47 & 3.03 & 2.34 & 2.85 & 3.35 \\
3 ($P_3/P_2$) & 3.01 & 3.01 & 3.01 & 3.01 & 3.01 & 3.09 & 3.30 & 3.43 & 3.51 & 3.52 \\
4 ($P_4/P_3$) & 4.00 & 4.00 & 4.00 & 4.00 & 4.00 & 4.11 & 4.16 & 4.17 & 4.27 & 4.43 \\
5 ($P_5/P_4$) & 5.00 & 5.01 & 5.00 & 5.02 & 5.29 & 5.18 & 5.24 & 5.01 & 5.31 & 5.29 \\
\hline
\end{tabular}        

\caption{Carreau law \eqref{eqn:carreau} with parameters $k_0 = 2$, $k_\infty = 0$, $\lambda = 2$ and $n = r \in \{ 1.2, 1.6, 2, 2.4, 2.8\}$. Experimental convergence rates for  different polynomial degrees $j$ (cf. theoretical estimate \eqref{eqn:err_uh} with $\theta=2$).}
        \label{table:conv-rate}
\end{table}

\section{Data-driven ICNN rheological models}
\label{sec:ICNN}
The aim of this section is two\-fold: (a) we briefly recall Input-Convex Neural Network (ICNN) scalar valued functions introduced in \cite{amos2017input} (cf. Section \ref{S:ICNN}); (b) we show how to employ ICNN to build  data-driven rheological models that are conformal with Assumptions ({\bf A}) (cf. Section \ref{S:data-driven_models}). The conformity with the assumptions is crucial to ensure the well-posedness of the associated differential problem governed by the obtained data-driven rheological model (see   \eqref{eqn:icnnstokes} below), thus building the numerical simulation upon solid mathematical foundations.  
%%%%
\subsection{Input-Convex Neural Networks}\label{S:ICNN}
We first consider a classical feedforward fully connected neural network $\text{NN}_{\bm{\theta}} (\mathbf{x}): \mathbbm{R}^d \rightarrow \mathbbm{R}$,  $d\geq 1$ where $\mathbf{x}=(x_1,\ldots,x_d)$ and the vector ${\bm{\theta}}$ contains the parameter identifying the neural network (weights and biases below).  We denote by $L$ the number of hidden layers.  We will make use of the following standard notations:\\
\begin{itemize}
    \item $w^{(l)}_{ij} \rightarrow$ weight connecting neuron $j$ in layer $l-1$ to neuron $i$ in layer $l$;
    \item $b^{(l)}_i \rightarrow$ bias of neuron $i$ in layer $l$;
    \item $h^{(l)}_i \rightarrow$ output of neuron $i$ in layer $l$, defined as: $h_i^{(l)} = \sigma_l ( \sum_j w_{ij}^{(l)}h_j^{(l-1)} + b_i^{(l)})$,  where $\sigma_l$ is the activation function at layer $l$.
\end{itemize}
The first and last layers correspond to the inputs and output of the neural network, respectively. Indeed, defining $y_j$ as the $j^{th}$ output, we have:

\begin{align*}
    &h_j^{(0)} = x_j \qquad (j = 1,..., d), \\
    &h_j^{(L+1)} = y:=\text{NN}_{\bm{\theta}}(\mathbf{x}).
\end{align*}

According to \cite{amos2017input},  the feedforward neural network function $\text{NN}_{\bm{\theta}}(\mathbf{x})$ turns into a {\it convex function} if the following conditions hold:
    \begin{itemize}
        \item[({\bf{C1}})] $w_{ij}^{(2:L+1)} \geq 0$;
         \item[({\bf{C2}})] $\sigma_l$ is convex and non-decreasing $\forall l \in \{1,\ldots,L+1\}$.
    \end{itemize}
In other words,  any feed-forward network function $\text{NN}_{\bm{\theta}} (\mathbf{x})$ can be re-worked into its convex counterpart, called $\ICNN (\mathbf{x})$, by choosing a non-decreasing (and convex) activation function and restricting its weights to be non-negative (for all but the first layer). 

Having in mind ({\bf{C1}}),  to guarantee the convexity of the function it becomes important to regulate the sign of each weight $w_{ij}^{(2:L+1)}$ during the training phase. If a weight $\overline{w} \in \{w_{ij}^{(k)}\}_{i,j,k}$, for $k= 1,\ldots,L+1$, becomes negative during the training process, it must be made positive before performing the successive training step. To reach this goal, various algorithms exist, such as nullifying negative weights by setting them to 0. However, motivated by the encouraging results presented in  \cite{sivaprasad2021curious}, in the sequel of the paper we employ the so-called {\em exponentiation algorithm}, described in Algorithm \ref{alg:exponentiation}. 
\begin{algorithm}[H]
    % \label{alg:exponentiation}
    \caption{Exponentiation}
    \begin{algorithmic}[1]\label{alg:exponentiation}
    
    \STATE INPUT: constant $\epsilon$
    \WHILE{training in progress}
        \STATE do training
        \FOR{$layer \in 2 : L+1$}
            \IF{$w < 0$}
                \STATE $w \leftarrow e^{w-\epsilon}$
            \ENDIF
        \ENDFOR
    \ENDWHILE

    \end{algorithmic}

\end{algorithm} 

The parameter $\epsilon$ entering in Algorithm \ref{alg:exponentiation} serves the purpose of constraining the updated weights, ensuring their proximity (from the right) to zero after exponentiation.
It should be noted that, as $\epsilon$ increases, there is a concurrent increase in the number of weights that are turned from negative to positive. 
For all the numerical experiments presented in this work, the value $\varepsilon = 30$ turned out to be a satisfactory compromise.\\

In what follows, with the aim of exploring the efficacy of Algorithm \ref{alg:exponentiation} in enforcing the convexity constraint, we present two simple examples. In the first case, we consider the approximation of the {\em non-convex} function $y(x) = |x| + \sin(x)$ by an Input-Convex Neural Network function $\ICNN$ and we compare it with the approximation obtained by employing a standard feedforward neural network function $\text{NN}_{\bm{\theta}}$. The training dataset is obtained by sampling the function $y$ in 100 randomly generated point belonging to the interval (-10,10). Both $\text{ICNN}_{\bm{\theta}}$ and $\text{NN}_{\bm{\theta}}$ were constructed with an architecture consisting of one input layer with 1 unit, followed by two hidden layers with 120 and 56 units, respectively, and finally one output layer with 1 unit. These networks were trained using the Mean Square Error (MSE) loss function, employing the Adam optimization algorithm over 3000 epochs.
The employed neural network architecture has proven to be the most effective one after undergoing numerous numerical tests.

\begin{figure}[H]
    \centering
    \includegraphics[width=0.65\linewidth]{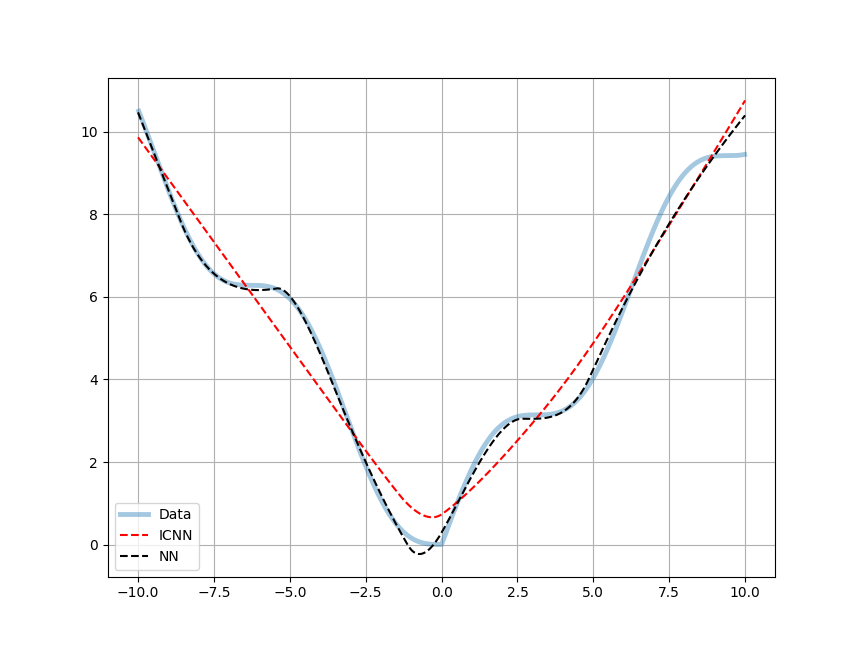}
    \caption{Approximation of $y(x) = |x| + \sin(x)$: comparison between $\ICNN$ and standard $\text{NN}_{\bm{\theta}}$ with same architecture, same training procedure.}
    \label{fig:icnn_vs_nn}
\end{figure}

In Figure \ref{fig:icnn_vs_nn}, a direct comparison between $\ICNN$ and $\text{NN}_{\bm{\theta}}$ is presented. We first observe that differently from $\text{NN}_{\bm{\theta}}$, the function $\ICNN$ is convex. It is also worth noting, for future use, that the convexity constraint plays the role of an implicit regularizing mechanism to prevent the overfitting of data.

\begin{figure}[H]
\centering
 \includegraphics[width=0.4\linewidth]{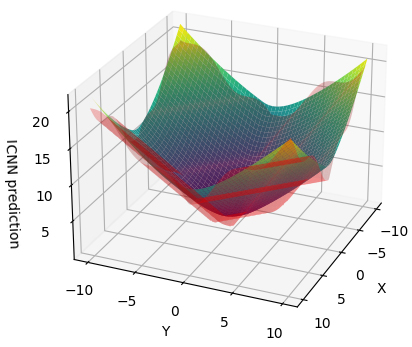}
 \qquad
 \includegraphics[width=0.4\linewidth]{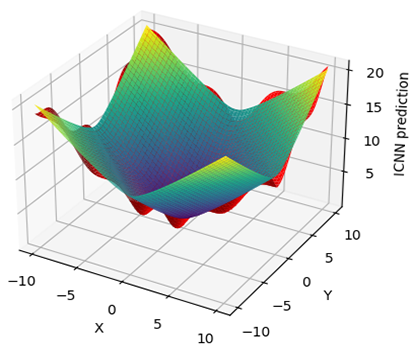}
    \caption{Approximation of the non-convex function $f(x,y)= |x| + |y| + \sin(x+y)$ (red colour) by $\ICNN$ (viridis): two different points of view.} \label{fig:2dnoinfo}
\end{figure}

As a second example, we consider $f(x,y)= |x| + |y| + \sin(x+y)$, a two-dimensional {\em non-convex} function which is approximated by an Input-Convex Neural Network function $\ICNN$ in the domain $\Omega=[-10,10]\times[-10,10]$. The architecture selected for $\ICNN$ comprised one input layer with 2 units, three hidden layers with respectively 120, 56 and 56 units, and one output layer with 1 unit. 
For the training, 1000 random points were randomly generated in $\Omega$. The network was trained using the Mean Square Error loss function and the Adam optimization algorithm over 10000 training epochs.
In Figure \ref{fig:2dnoinfo}, two different perspectives of the approximation of the non-convex function $f$ by ICNN are presented. We observe that the ICNN effectively captures the global trend of the function while disregarding, because of the convexity constraint, the contribution from the oscillating term $\sin(x+y)$.

At last, before detailing the data-driven construction of rheological models, let us remark that in case one needs to select among {\em convex} or {\em concave} approximation to best fit a given dataset, the simple procedure presented in Algorithm \ref{alg:convexorconcave} can be employed. The idea is straightforward and it is based on employing the given dataset, say $\{(x_i,y_i)\}_{i=1}^N$ and its modified version $\{(x_i,-y_i)\}_{i=1}^N$ to train two Input-Convex Neural Networks, say $\ICNN^+$ and $\ICNN^-$.
Finally, comparing the values of the loss functions associated to $\ICNN^+$ and $\ICNN^-$ yields the best fit. Note that while $\ICNN^-$ is a convex approximation of the modified data set $\{(x_i,-y_i)\}_{i=1}^N$, the opposite function $-\ICNN^-$ is a {\em concave} approximation of the {\em original} dataset $\{(x_i,y_i)\}_{i=1}^N$.

\begin{algorithm}[H]
    \caption{Convex/concave selection}
    \label{alg:convexorconcave}
    \begin{algorithmic}
    \STATE INPUT: dataset $\{(x_i,y_i)\}_{i=1}^N$
    \WHILE{training in progress}
        \STATE train $\ICNN^+(\mathbf{x})$ with data $\{(x_i,y_i)\}_{i=1}^N$ and $\ICNN^-(\mathbf{x})$ with data $\{(x_i,-y_i)\}_{i=1}^N$
    \ENDWHILE
    \IF {loss($\ICNN^+(\mathbf{x})$) <
        loss($\ICNN^-(\mathbf{x})$)}
        \STATE $\ICNN=\ICNN^+(\mathbf{x})$
    \ELSE
        \STATE $\ICNN=-\ICNN^-(\mathbf{x})$
    \ENDIF
    \end{algorithmic}
\end{algorithm} 
\subsection{Data-driven rheological models: the use of ICNNs}\label{S:data-driven_models}
We are now ready to discuss how to employ convex neural networks to build data-driven rheological models of generalized non-Newtonian fluids, whose viscous stress, we recall, is given by 
\begin{equation}\label{eq:tau} 
\tau(\varepsilon(\mathbf{u})) = k( |\varepsilon(\mathbf{u}) |)\varepsilon(\mathbf{u}),
\end{equation}
where the function $k$ is convex or concave depending on the nature of the fluid (shear-thinning or shear-thickening). 

More precisely, given experimental measurements, say $\{(x_i,y_i)\}_{i=1}^N$ with $x=|\varepsilon(\mathbf{u})|$ and $y=k(|\varepsilon(\mathbf{u})|)$, we employ  Algorithm \ref{alg:convexorconcave} that  automatically returns the best {\em convex/concave} fit $\ICNN$ of the viscosity $k$, thus providing the following {\em ICNN stress tensor}
\begin{equation}\label{eq:tau_theta}
\tau_{\bm{\theta}}(\varepsilon(\mathbf{u})) = \ICNN( |\varepsilon(\mathbf{u}) |)\varepsilon(\mathbf{u}).
\end{equation}
In the sequel, we will refer to the function $\ICNN$ as the {\em ICNN viscosity}.
 
To test this procedure on real rheological data, we consider now experimental rheological measurements for a set of  shear-thinning aqueous solutions of Xanthan gum with sodium chloride addition, available in the dataset \cite{mrokowska2023dataset}.
The solutions differs by molar concentration $M$ of NaCl ranging between 0 and 0.7, while the Xanthan gum content was the same for all solutions (1 g/L). The apparent viscosity for the different solutions has been experimentally characterized using a rotational rheometer (see \cite{mrokowska2023dataset} for details).   
For each experimental dataset, we adopt our strategy to learn the corresponding data-driven ICNN rheological model.

\begin{figure}[H]
    \centering
    \includegraphics[width=0.989\linewidth]{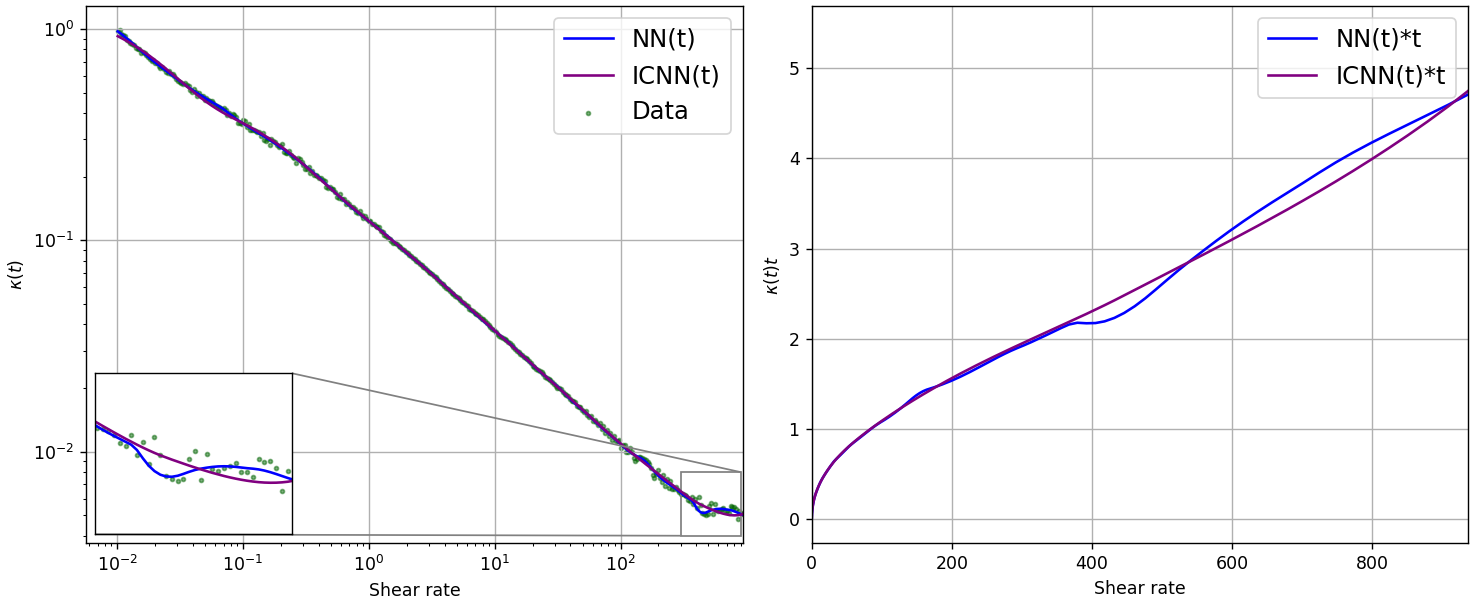}
    \caption{Comparison between $\ICNN$ and a standard feed-forward neural network $NN_{\bm{\theta}}$ in the approximation of real rheological measures: viscosity curves $k(t)$ (left) and function $k(t)t$ (right), with $k\in \{\ICNN, \text{NN}_{\bm{\theta}}\}$. DATASET: NaCL05\_XG}
    \label{fig:comparison_NNvsICNN}
\end{figure}

We first consider the solution with a NaCl molar concentration $M=0.5$ (labelled  in the sequel as NaCL\_05+XG) and we compare the learned rheological model obtained by employing  a standard feed-forward neural network $\text{NN}_{\bm{\theta}}$ with the one obtained with convex neural networks $\ICNN$. An identical architecture is used for both neural networks: specifically, a configuration of $1 \times 120 \times 56 \times 1$. Consistency in training was maintained by using the same dataset and MSE loss function for both networks, with the training process spanning over  $N=100000$ epochs.
In Figure \ref{fig:comparison_NNvsICNN} (left), the resulting viscosities corresponding to the two neural networks are reported, together with the dataset they are trained on. The ability of $\ICNN$ to avoid overfitting can be clearly appreciated. Moreover, from Figure \ref{fig:comparison_NNvsICNN} (right) the lack of monotonicity of $\text{NN}_{\bm{\theta}}(t)t$ for the standard feed-forward neural network clearly implies the violation of  Assumption $\textbf{(A)}_1$, while $\ICNN(t)t$ remains monotone (see also  Algorithm \ref{alg:ass_verifier} below for a systematic validation of Assumptions $\textbf{(A)}$).

\begin{figure}[H]
    \centering
    \subfloat[NaCL\_00+XG]{
    \includegraphics[width=0.4\linewidth]{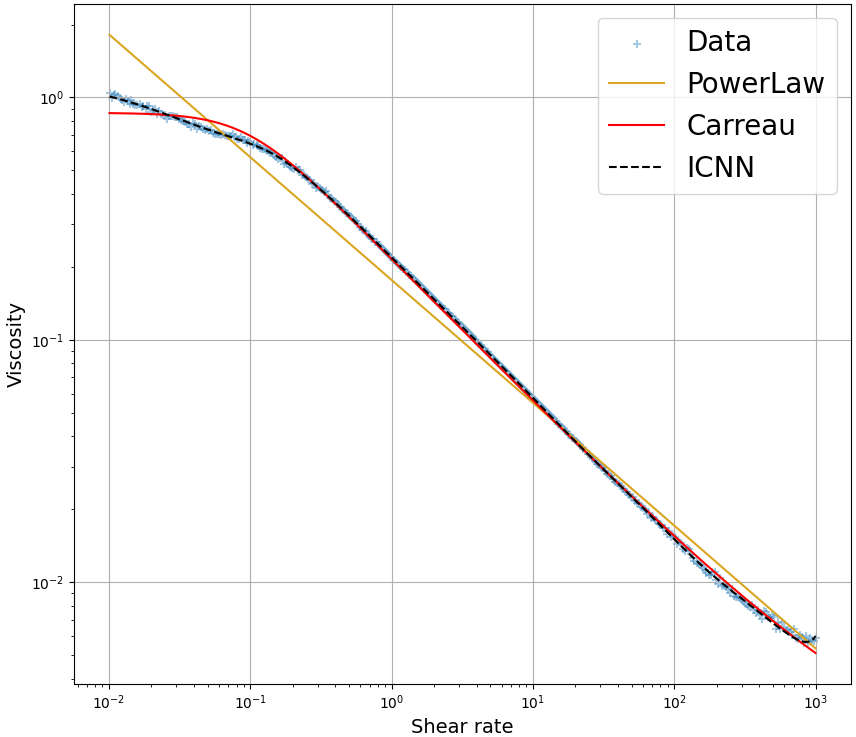}
    }
    \subfloat[NaCL\_01+XG]{
    \includegraphics[width=0.4\linewidth]{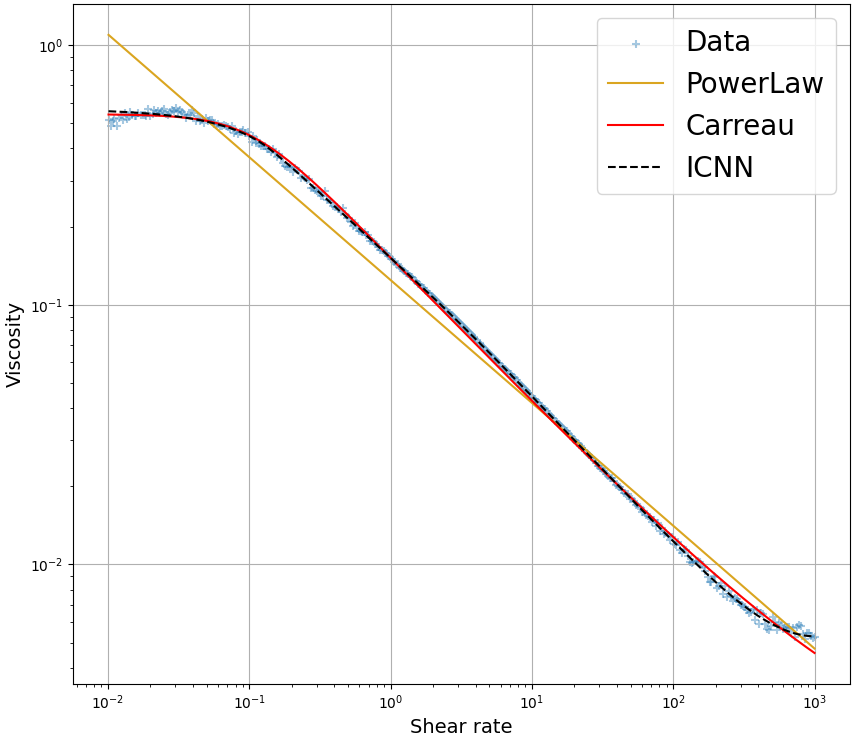}
    }
    \\
    \subfloat[NaCL\_05+XG]{
    \includegraphics[width=0.4\linewidth]{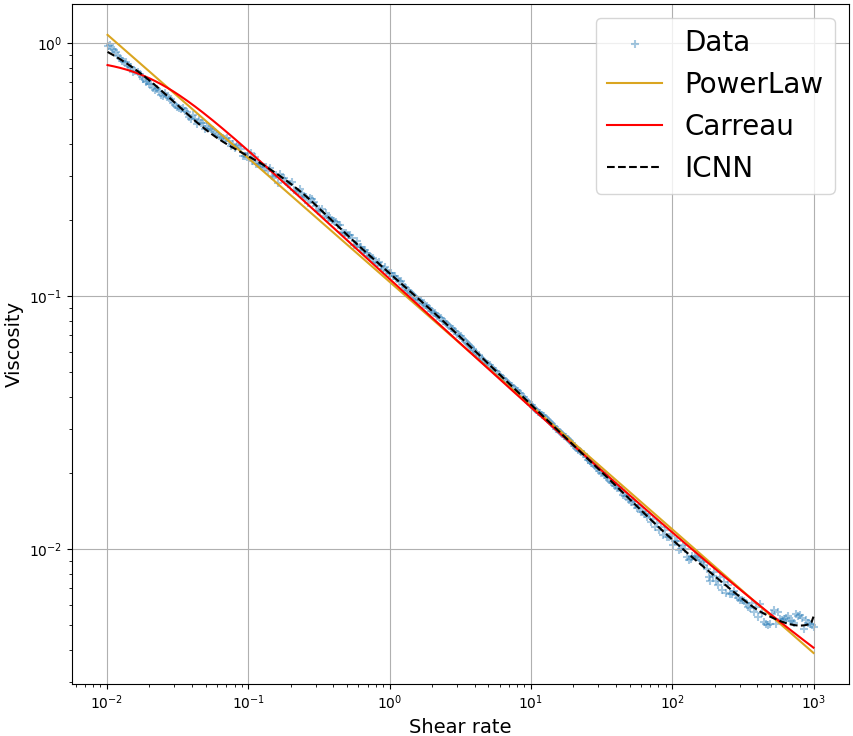}
    }
    \subfloat[NaCL\_07+XG]{
    \includegraphics[width=0.4\linewidth]{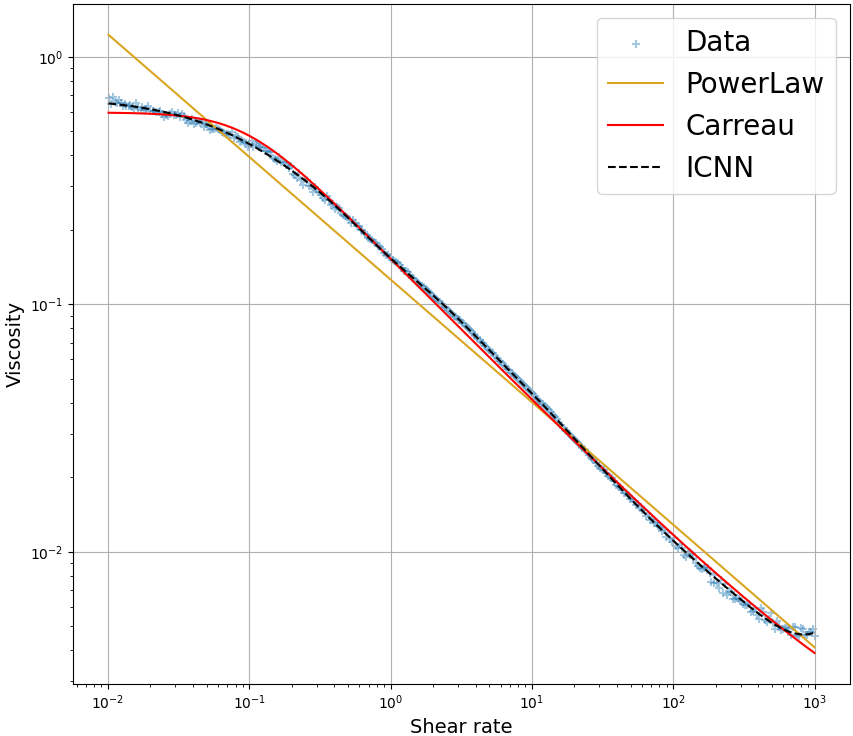}
    }
    \caption{Comparison between our data-driven ICNN procedure, Carreau and Power law  rheological models on 4 shear-thinning aqueous solutions of Xanthan gum with sodium chloride addition \cite{mrokowska2023dataset} . }
    \label{fig:ICNN_solutions}
\end{figure}

The proposed data-driven strategy is then employed to discover the rheology of the aqueous solution with NaCl molar concentration $M=0, 0.1, 0.5, 0.7$. The results are presented in Figure \ref{fig:ICNN_solutions}, where the ICNN viscosities are compared with the ones obtained employing  a-priori chosen parametric models, namely Carreau and Power law. The results shows the capability of  ICNN to automatically learn most of the trends presented in the datasets, including non-standard behaviour such as the one displayed by the NaCL\_00+XG dataset (see Figure \ref{fig:ICNN_solutions} (a)) for low shear rate values, that cannot be captured by the a-priori functional dependence prescribed by standard rheological models. 
A quantitative comparison in terms of both the Root Mean Squared Error (RMSE) and the  coefficient of determination $R^2$ is presented in Table \ref{table:comparison}. These results confirm the superior performances of the proposed  data driven ICNN rheological model.

\begin{table}[H]
    \centering
    \begin{tabular}{|c|cc|cc|cc|}
        \hline 
        & \multicolumn{2}{c|}{Power Law} & \multicolumn{2}{c|}{Carreau} & \multicolumn{2}{c|}{ICNN} \\ 
        \hline 
        Case  & RMSE & $R^2$ & RMSE & $R^2$ & RMSE & $R^2$ \\
        \hline
        NaCL\_00+XG      & 0.027362 & 0.72873 & 0.001160 & 0.98850 & 0.000037 & 0.99963 \\
        NaCL\_01+XG      & 0.012810 & 0.67595 & 0.000128 & 0.99761 & 0.000099 & 0.99750 \\
        NaCL\_05+XG      & 0.000859 & 0.98476 & 0.000731 & 0.98703 & 0.000035 & 0.99938 \\
        NaCL\_07+XG   & 0.012906 & 0.72695 & 0.000310 & 0.99344 & 0.000046 & 0.99903 \\
        \hline
    \end{tabular}
    \caption{Fitting comparison between Power Law, Carreau, and ICNN models.}
    \label{table:comparison}
\end{table}

Finally, we discuss how to verify that the learned ICNN viscosity satisfies Assumptions \textbf{(A)}. This is crucial  for the existence and uniqueness of the solution of the associated ICNN non-Newtonian Stokes (see \eqref{eqn:icnnstokes}). To reach this goal, we propose the use of {\bf Algorithm \ref{alg:ass_verifier}} which numerically checks  whether the trained $\ICNN$ is in compliance 
with the Assumptions \textbf{(A)} by fine-tuning (through a minimization process) the values of constants $r \in (1, \infty)$, $\alpha \in [0,1]$, and $C, M > 0$. 
%in 
%such a way that $C$ and $\alpha$ are minimized, while $M$ is maximized.

\begin{algorithm}[H]
    \caption{(Check the validity of Assumptions  \textbf{(A)}) }
    \begin{algorithmic}[1]\label{alg:ass_verifier}
    \STATE initial guess $C = C_0$, $\alpha = \alpha_0$, $r = r_0$, $M=M_0$; $t_1,...,t_N \in (0,\infty)$, $s_1,...,s_N \in (0, \infty)$
    \STATE $k=0$, max iterations $K > 0$
    \WHILE{ $k \leq K$}
        \STATE Define $f_1 = \sum_{j=1}^{N} C_k[ t_j^{\alpha_k}(1+t_j)^{1-\alpha_k}    ]^{r_k-2} - \text{ICNN}_{\bm{\theta}}(t_j) = \sum_{j=1}^{N} f_{1,j}$
        
        \STATE Define $f_2 = \sum_{j=1}^{N}\sum_{i=1}^{N} C_k | t_j-s_i | [(t_j+s_i)^{\alpha_{k}}(1+t_j+s_i)^{1-\alpha_{k}}]^{r_{k}-2} - | \text{ICNN}_{\bm{\theta}}(t_j)t_j- \text{ICNN}_{\bm{\theta}}(s_i)s_i | = \sum_{j=1}^{N}\sum_{i=1}^{N} f_{2,ij}$

        \STATE Define $f_3 = \sum_{j=1}^{N}\sum_{i=1}^{N}\text{ICNN}_{\bm{\theta}}(t_j)t_j- \text{ICNN}_{\bm{\theta}}(s_i)s_i - M_k(t_j-s_i)[(t_j+s_i)^{\alpha_{k}}(1+t_j+s_i)^{1-\alpha_{k}}]^{r_{k}-2} = \sum_{j=1}^{N}\sum_{i=1}^{N} f_{3,ij}$
        \STATE Compute $F = f_1 + f_2 \mathds{1}_{\{|s/t - 1| \leq 1\}} + f_3 \mathds{1}_{\{t \geq s\}} = \sum_{j=1}^{N} f_{1,j} + \sum_{j=1}^{N}\sum_{i=1}^{N} f_{2,ij} \mathds{1}_{\{|s_i/t_j - 1| \leq 1\}} + \sum_{j=1}^{N}\sum_{i=1}^{N} f_{3,ij} \mathds{1}_{\{t_j \geq s_i\}} $
        \STATE Minimize $F$ with respect to $C_k,M_k, \alpha_k,r_k$ imposing $f_{1,j} \geq 0 \quad \forall t_j$, and $ f_{2,ij} \geq 0 \quad \forall t_j, s_i$ such that $|s_i/t_j - 1| \leq 1$, and $ f_{3,ij} \geq 0 \quad \forall t_j, s_i$ such that $t_j \geq s_i$
        \STATE $k = k+1$
    \ENDWHILE
    \end{algorithmic}
\end{algorithm}

In Table \ref{table:constants_real} we collect the output of Algorithm \ref{alg:ass_verifier} applied to the ICNN approximation of the dataset in \cite{mrokowska2023dataset}. The initial values of the parameters $C_0$, $\alpha_0$, $r_0$, and $M_0$ were set to 10, 0.5, 1.5, and 10, respectively. It is worth noting that the specific choice of these initial values did not significantly impact the performance of the minimization algorithm.
The minimization problem in Step 8 of Algorithm \ref{alg:ass_verifier} is solved using the \texttt{best1bin} algorithm implemented in the \texttt{differential\_evolution} function of the \texttt{scipy.optimize} library. The \texttt{differential\_evolution} function was favored over traditional optimization algorithms due to its multi-process capabilities, enhancing computational efficiency, and its proficiency in avoiding local minima, ensuring more robust solutions. The arrays $\{t_i\}$ and $\{s_j\}$,  with $i, j=1, \dots,100$ were generated uniformly within the range of 0 to the maximum shear rate value stored in the dataset. The value $K$ was selected to be arbitrarily high (over 10,000 epochs) to ensure robust performance during the iteration process.

\begin{table}[H]
    \centering 
        \begin{tabular}{|c |c |c |c |c|}
        \hline
     
         & NaCL\_00+XG &  NaCL\_01+XG & NaCL\_05+XG & NaCL\_07+XG \\
        \hline
        $C$ & 1.002 & 0.538 & 1.052 & 0.635\\
        $\alpha$ & 0.015 & 0.010  & 0.013  & 0.015\\
        $r$ & 1.628 & 1.452 & 1.453 &  1.624\\
        $M$ & 0.050 & 0.098 & 0.053 & 0.038\\
        \hline
        \end{tabular}
        \\[10pt]
        \caption{Optimized constants  ensuring that the \text{ICNN} viscosities trained on the datasets \cite{mrokowska2023dataset} satisfy Assumptions \textbf{(A)}.}
        \label{table:constants_real}
\end{table}

\section{Non-Newtonian Stokes equations with neural networks}
\label{sec:ICNNStokes}

{
In this section, we study the finite element approximation of  the non-Newtonian Stokes equations \eqref{eqn:weak_stokes}, where the viscous stress $\tau$ is replaced by its ICNN approximation $\tau_\theta$ (cf. \eqref{eq:tau} and \eqref{eq:tau_theta}). The corresponding problem is named {\em ICNN non-Newtonian Stokes problem} and its  
weak formulation reads as follows: find $(\mathbf{u},p) \in \mathbf{X} \times M$ such that:
\begin{equation}\label{eqn:icnnstokes}
    \begin{cases}
        \begin{aligned}
            &\textstyle \int_\Omega \text{ICNN}_{\bm{\theta}}(\left\vert \varepsilon(\mathbf{u})\right\vert)\varepsilon(\mathbf{u}):\varepsilon(\mathbf{w})  - \int_\Omega p\nabla \cdot \mathbf{w} = \int_\Omega \mathbf{f} \cdot \mathbf{w}\quad \forall \mathbf{w}\in \mathbf{X}
        \end{aligned} \\
        \textstyle \int_\Omega q \nabla \cdot \mathbf{u} = 0\quad \forall q\in M.
    \end{cases}
\end{equation}
In particular, we are interested  in studying the  convergence properties of the finite-element approximation of \eqref{eqn:icnnstokes} towards the solution of \eqref{eqn:weak_stokes}. To this aim, it will be crucial to incorporate the contribution from the approximation $\tau_{\bm{\theta}}$ of the stress tensor $\tau$ (see Theorem \ref{prop:ourprop} below) and combine it with the results of Theorem \ref{thm:convergence} (see the inequality  \eqref{eq:fem+ICNN} below). 

To ease the presentation,  we train the data-driven ICNN rheological model on datasets obtained by sampling the Carreau law \eqref{eqn:carreau} to obtain 100 shear-rate values randomly generated over the interval (0,70). More precisely, 
to explore a variety of different scenarios, multiple datasets based on \eqref{eqn:carreau} have been produced, choosing $k_\infty=0$, $k_0=2$, $\lambda=2$, and different values of the exponent $n=r=1.2, 1.6, 2, 2.4, 2.8$. The range of $n$ allows to consider both pseudo-plastic and dilatant regimes.
Since, depending on the regime, the resulting viscosity $k$ may be either convex (for $1 < r \le 2$) or concave (for $r \ge 2$), we employed Algorithm \ref{alg:convexorconcave} to automatically select the best convex/concave fit. 
The architecture of the employed ICNN is $1\times120\times56\times1$ and 20000 epochs of the Adam optimization algorithm have been run to minimize the loss function given by the Mean Square Error. 
The resulting $\ICNN$ functions approximating the viscosity $k$ for  different values of $n$ are reported in Figure \ref{fig:icn}, while Table \ref{table:icnn-errors} collects the associated $L^2$ error.

\begin{figure}[H]
    \centering
    \includegraphics[width=0.6\linewidth]{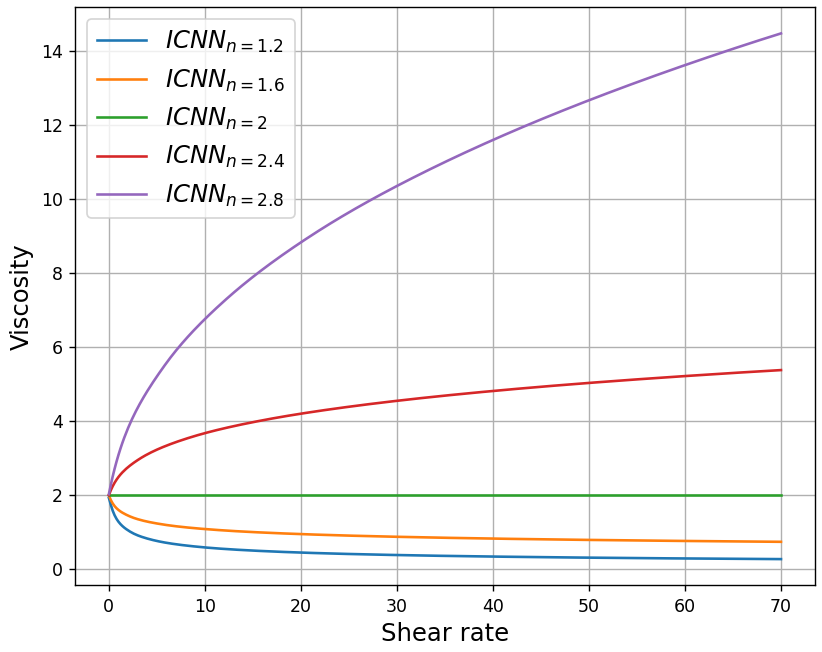}\\
    \caption{$\ICNN$ approximations of the  Carreau law \eqref{eqn:carreau} with parameters: $k_\infty=0,k_0=2,\lambda=2,n=1.2, 1.6, 2, 2.4, 2.8$.}
    \label{fig:icn}
\end{figure}

\begin{table}[H]
    \centering 
        \begin{tabular}{|c | c  |}
        \hline
         &  $L^2$-error  \\
        \hline
        $\text{ICNN}_{n=1.2}$ & 4.1e-3 \\
        $\text{ICNN}_{n=1.6}$ & 9.6e-4 \\
        $\text{ICNN}_{n=2.0}$ & 1.1e-6  \\
        $\text{ICNN}_{n=2.4}$ & 6.2e-4 \\
        $\text{ICNN}_{n=2.8}$ & 6.9e-3 \\
        \hline
        \end{tabular}
        \\[10pt]
        \caption{$L^2$-error between the Carreau viscosity $k$ and the ICNN viscosity for different values of $n$, measured on the interval $(0,70)$. }
        \label{table:icnn-errors}
\end{table}

Finally, in Table \ref{table:constants} we report the values of the optimized constants, obtained by Algorithm \ref{alg:ass_verifier}, to ensure the validity of Assumptions ({\bf{A}}) and thus the well-posedness of the ICNN non-Newtonian Stokes problem \eqref{eqn:icnnstokes}.
In particular, we note that the obtained values of $r$ are closed, as expected, to the corresponding values of $n$. 
\begin{table}[H]
    \centering 
        \begin{tabular}{|c| c| c| c| c| c|}
        \hline
     
         & $\text{ICNN}_{n=1.2}$ &  $\text{ICNN}_{n=1.6}$ &  $\text{ICNN}_{n=2.0}$ & $\text{ICNN}_{n=2.4}$ & $\text{ICNN}_{n=2.8}$ \\
        \hline
        $C$      & 1.600 & 1.786 & 2.000   & 2.379  & 2.753 \\
        $\alpha$ & 0.099 & 0.118 & 0.022   & 0.182  & 0.162 \\
        $r$      & 1.176 & 1.589 & 2.000   &  2.403 & 2.807 \\
        $M$      & 0.859 & 0.913 & 1.000   & 0.369  &  0.136 \\
        \hline
        \end{tabular}
        \\[10pt]
        \caption{Optimized constants obtained by Algorithm \ref{alg:ass_verifier}: ICNN functions  satisfy Assumptions {\textbf{(A)}}.} 
        \label{table:constants}
\end{table}

\subsection{Convergence results}
Let ($\overline{\mathbf{u}},
\overline{p}$) be the solution of the non-Newtonian Stokes problem \eqref{eqn:weak_stokes} where the viscosity is governed by the Carreau law \eqref{eqn:carreau} with parameters: $k_\infty=0,k_0=2,\lambda=2$ and $n\in \{1.2, 1.6, 2, 2.4, 2.8\}$.  
Let us denote by ($\mathbf{u}_{h,ICNN}$, $p_{h,ICNN}$) the finite element approximation of  ICNN-Stokes equations \eqref{eqn:icnnstokes} with ICNN viscosity obtained in the previous section.
In Figure \ref{fig:nn_errors}, we plot the velocity error  $\|\overline{\mathbf{u}} - \mathbf{u}_{h,ICNN} \|_{[W^{1,r}(\Omega)]^2}$  and the  pressure error $\|\overline{p} - p_{h,ICNN} \|_{L^{r'}(\Omega)}$ as the mesh is refined, for $r \in \{ 1.2,1.6,2, 2.4, 2.8\}$ (recall $r=n $). These results have been obtained with Taylor-Hood finite-elements with $j=2$.

 The plots in Figure \ref{fig:nn_errors} reveal two key findings. First, the errors 
steadily decrease to a minimum in all different scenarios (apart the trivial case $r=2$). Secondly, the value of the plateau is related to the approximation error of $\text{ICNN}$ when trained to fit the Carreau law (see Table \ref{table:icnn-errors}). 
More precisely, denoting $\varepsilon_{\text{ICNN},u}$ and $\varepsilon_{\text{ICNN},p}$ as
$$
 \varepsilon_{\text{ICNN},u}:=\lim_{h \rightarrow 0} \| \overline{\mathbf{u}} - \mathbf{u}_{h,ICNN} \|_{[W^{1,r}(\Omega)]^2} 
$$
$$
\varepsilon_{\text{ICNN},p}:=\lim_{h \rightarrow 0} \| \overline{p} - p_{h,ICNN} \|_{L^{r'}(\Omega)} 
$$
we note a correlation between 
$\varepsilon_{\text{ICNN},u}$, $\varepsilon_{\text{ICNN},p}$ and the errors reported in Table \ref{table:icnn-errors}, as it is evident from Figure \ref{fig:vareps}, where we observe that an improvement of the accuracy of $\text{ICNN}$  corresponds to lower values of $\varepsilon_{\text{ICNN},u}$ and $\varepsilon_{\text{ICNN},p}$.

\begin{figure}[H]
    \centering
    \subfloat[Case \( n=r=1.2 \)]{
        \includegraphics[width=0.75\linewidth]{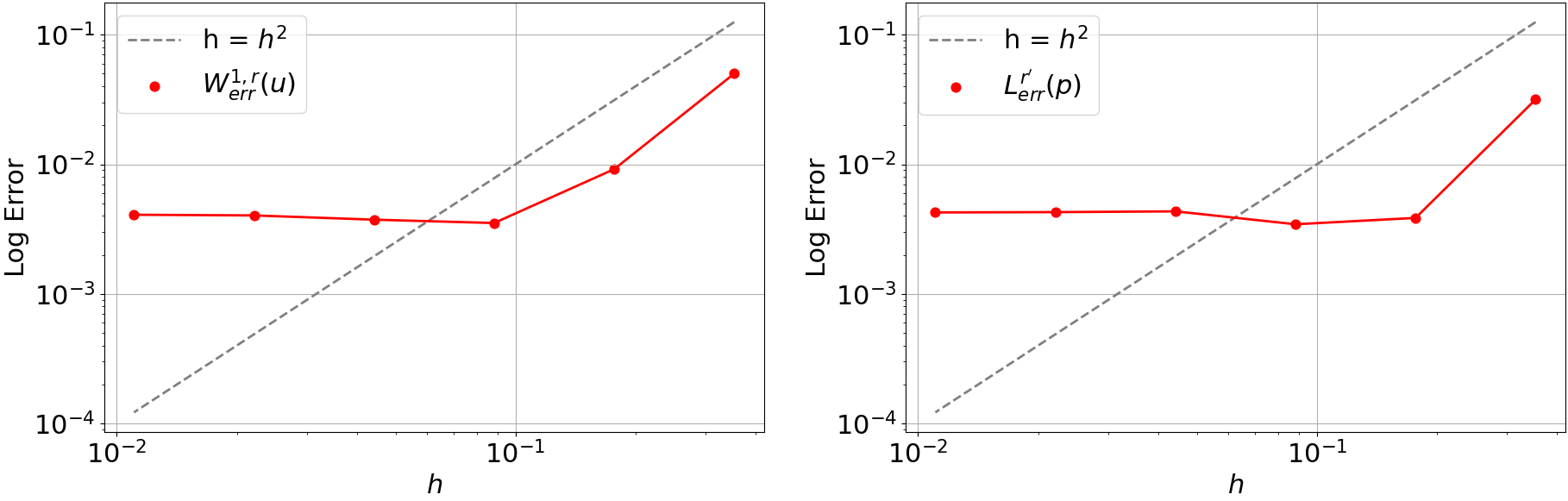}
    }
    \\
    \subfloat[Case \( n=r=1.6 \)]{
        \includegraphics[width=0.75\linewidth]{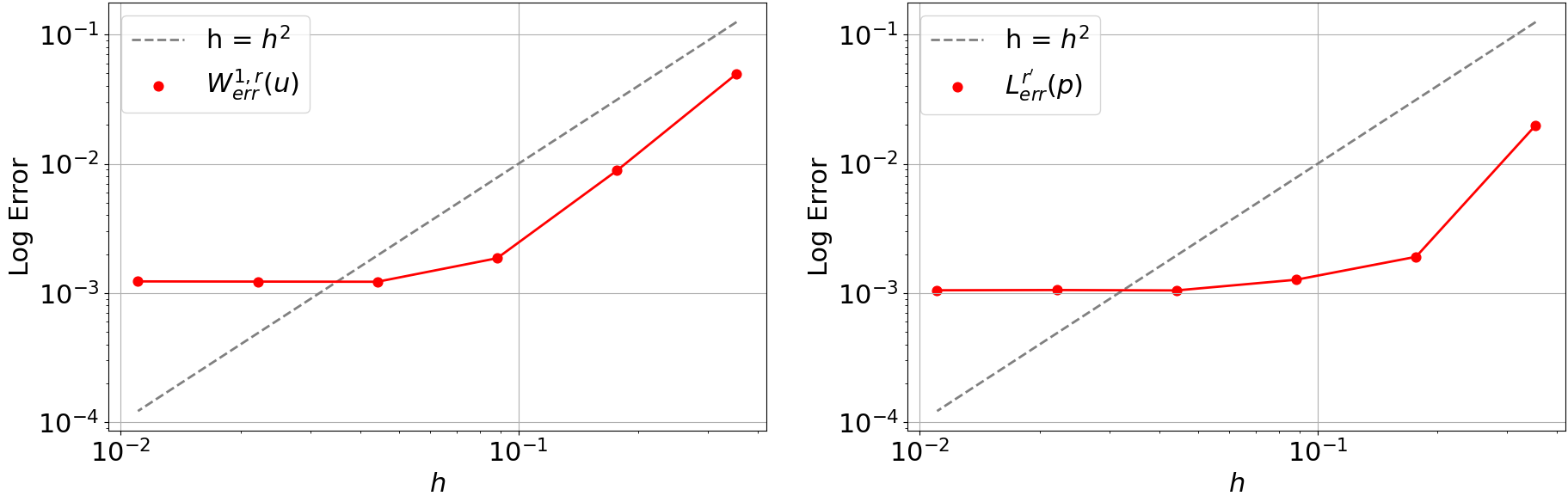}
    }
    \\ % Added line break for better formatting
    \subfloat[Case \( n=r=2 \)]{
        \includegraphics[width=0.75\linewidth]{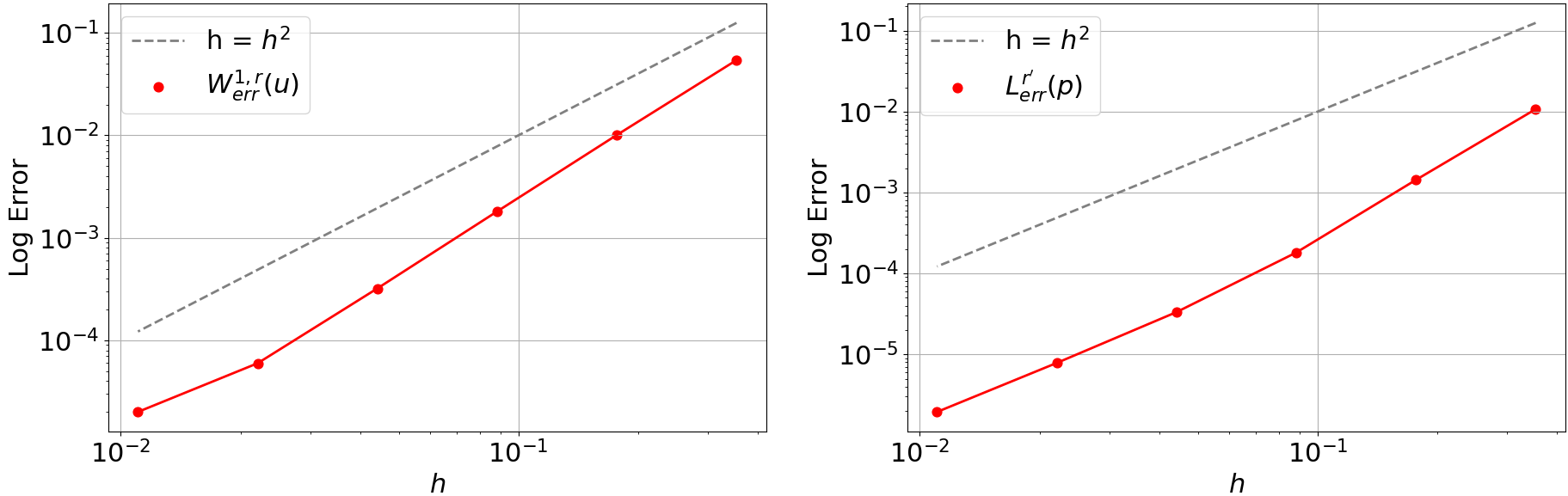}
    }
    \\
    \subfloat[Case \( n=r=2.4 \)]{
        \includegraphics[width=0.75\linewidth]{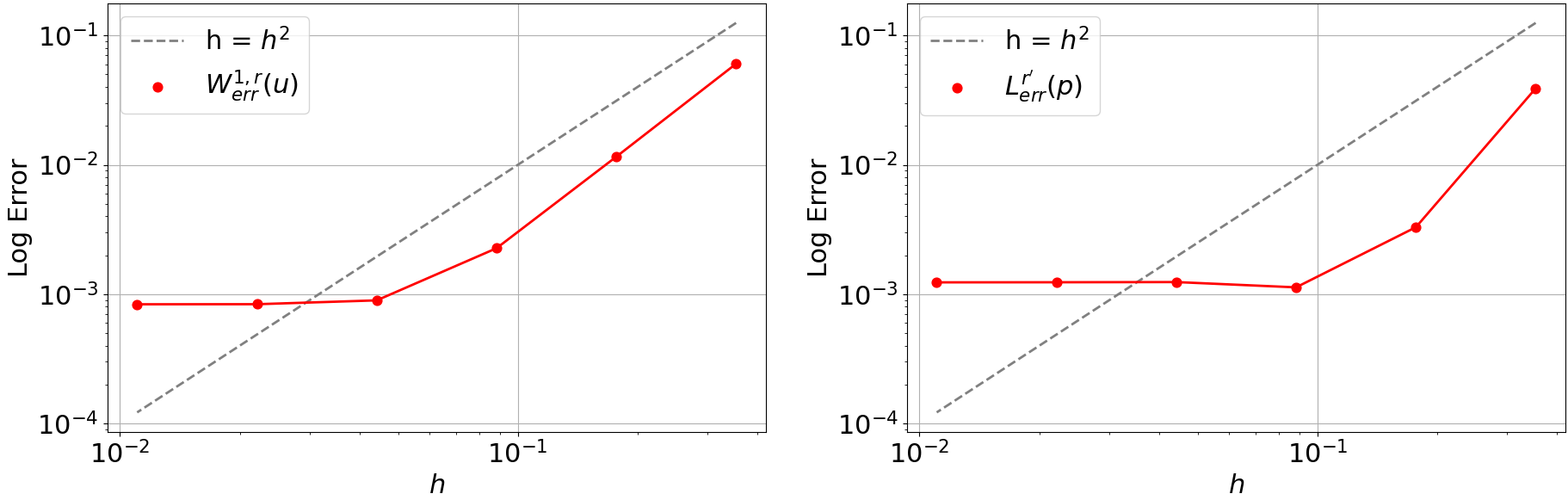}
    }
    \\
    \subfloat[Case \( n=r=2.8 \)]{
        \includegraphics[width=0.75\linewidth]{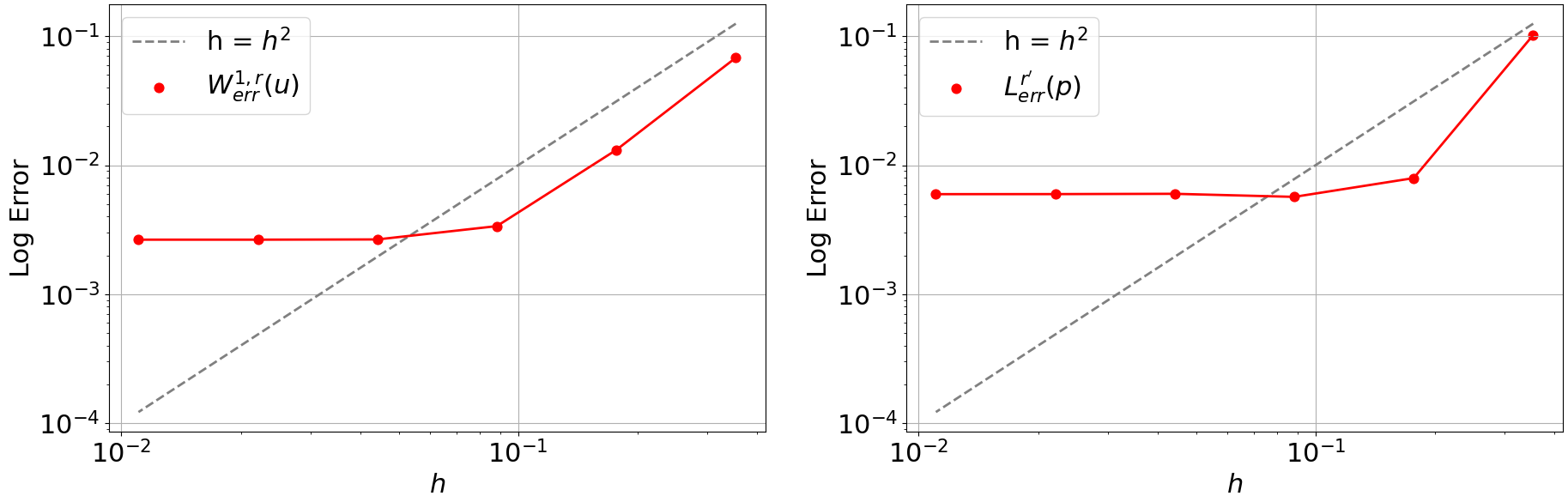}
    }
    \caption{Velocity error \(\|\overline{\mathbf{u}} - \mathbf{u}_{h,ICNN} \|_{[W^{1,r}(\Omega)]^2}\) (left) and pressure error \(\|\overline{p} - p_{h,ICNN} \|_{L^{r'}(\Omega)}\)for \( r \in \{ 1.2,1.6,2, 2.4, 2.8\}  \).}
    \label{fig:nn_errors}
\end{figure}

\begin{figure}[H]
    \centering
    \includegraphics[width=0.9\linewidth]{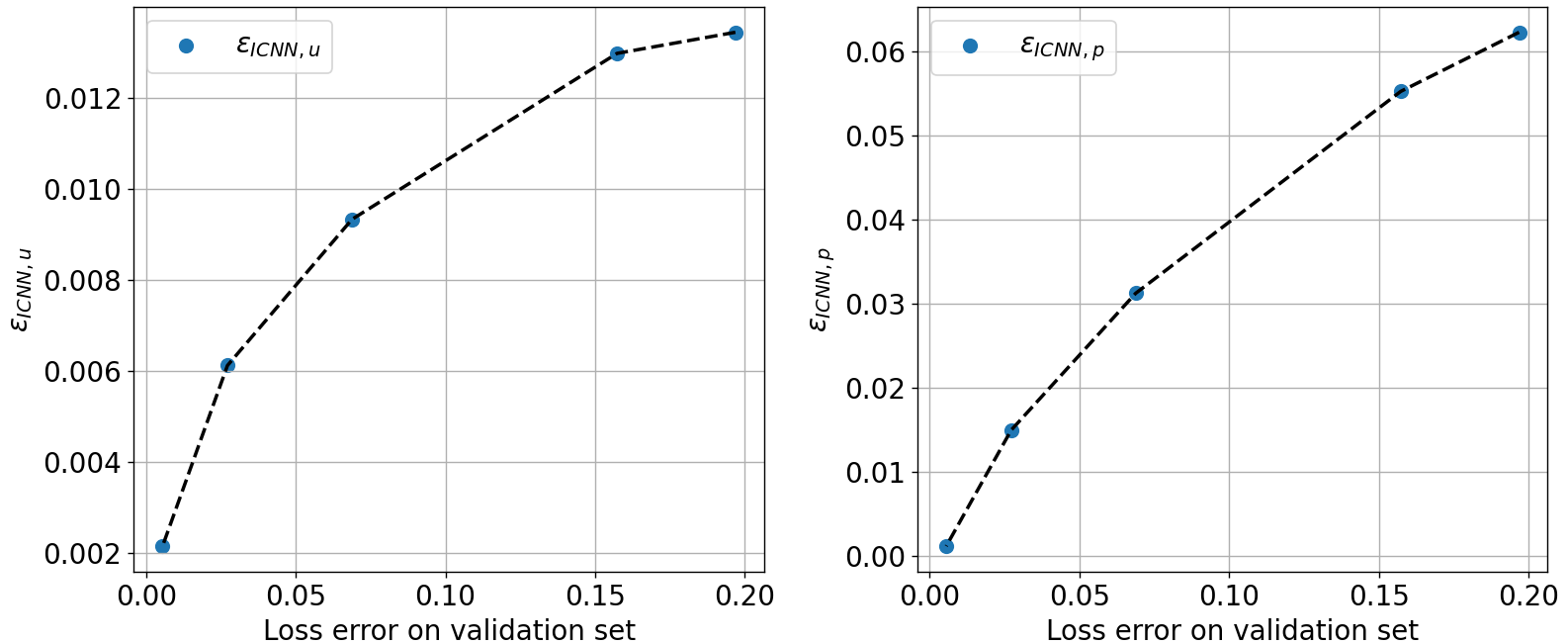}
    \caption[Behaviour of $\varepsilon_{\text{ICNN},u}$ and $\varepsilon_{\text{ICNN},p}$]{$\varepsilon_{\text{ICNN},u}$ and $\varepsilon_{\text{ICNN},p}$ against  the accuracy of $\text{ICNN}$.}
    \label[type]{fig:vareps}
\end{figure}

\begin{remark} As mentioned earlier, for $r=2$ no plateau is observed and asymptotic optimal convergence is achieved. This is due to the fact that, in this case, the function $k(\cdot)$ reduces to a constant and can be approximated exactly by $\text{ICNN}$.
\end{remark}

The following perturbation result paves the way for a deeper  understanding of the above numerical findings. 
 
\begin{theorem}[Perturbation result]\label{prop:ourprop}
    Let $\mathbf{X}_l = [W^{1,r_l}(\Omega)]^2$, $M_l=L^{r_l'}(\Omega)$, $l=1,2$. Let $A_l: \mathbf{X}_l \rightarrow \mathbf{X}_l^*$, $l=1,2$ be defined as:
    $
        \left\langle A_l\mathbf{v}, \mathbf{w} \right\rangle_{\mathbf{X}^*_l} := \int_\Omega k_l(\left\vert \varepsilon(\mathbf{v})\right\vert)\varepsilon(\mathbf{v}):\varepsilon(\mathbf{w}),
    $ for any $\mathbf v,\mathbf w\in \mathbf{X}_l$, where $k_l$ satisfies Assumptions {\bf (A)} for $l=1,2$. 
    Let $(\mathbf{u}_l, p_l)$ be the solution to:
\begin{equation*}
    \begin{cases}
        \left\langle A_l\mathbf{u}_l, \mathbf{w}\right\rangle_{\mathbf{X}_l^*} - \left\langle p_l, \nabla \cdot \mathbf{w}\right\rangle  = 
        \left\langle \mathbf{f},\mathbf{w}\right\rangle  &\qquad \forall \mathbf{w} \in \mathbf{X}_l,\\
        \left\langle \nabla \cdot \mathbf{u}_l, q   \right\rangle  = 0 &\qquad \forall q \in M_l,\\
    \end{cases}
\end{equation*}
for $l = 1,2$. Set $r:=\min_{l=1,2}\{r_l\}$, denote by $j$ the index such that $r_j=\max_{l=1,2}\{r_l\}$ and by $i$ the other index. 
If there exists $q\in[r',\infty]$ such that $\Vert \varepsilon(\mathbf{u}_j)\Vert_{[L^{s}(\Omega)]^4}<+\infty$, where $s=\frac{1}{1-\frac 1 q-\frac 1r}$ if $q>r'$, $s=\infty$ otherwise, and $\| k_1(| \varepsilon(\mathbf{u}_j) |) - k_2(| \varepsilon(\mathbf{u}_j) |)\|_{L^q(\Omega)}<+\infty$, then the following inequalities hold:
\begin{eqnarray}
  &  \| \mathbf{u}_1 - \mathbf{u}_2 \|_{[W^{1,r}(\Omega)]^2} \leq & C \| k_1(| \varepsilon(\mathbf{u}_j) |) - k_2(| \varepsilon(\mathbf{u}_j) |) \|_{L^q(\Omega)}\| \varepsilon(\mathbf{u}_j) \|_{[L^{s}(\Omega)]^4},\label{eqn:result}
  \end{eqnarray}
  and, when $r\leq 2$,
  \begin{eqnarray}
     \| p_1 - p_2 \|_{L^{r'}(\Omega)} \nonumber&\leq C 
     \| k_1(| \varepsilon(\mathbf{u}_j) |) - k_2(| \varepsilon(\mathbf{u}_j) |)\|_{L^q(\Omega)}  
    \| \varepsilon(\mathbf{u}_j) \|_{[L^{s}(\Omega)]^4}\\ &\qquad +\ C_i 
     \| k_1(| \varepsilon(\mathbf{u}_j) |) - k_2(| \varepsilon(\mathbf{u}_j) |)\|_{L^q(\Omega)}^{r-1}  
    \| \varepsilon(\mathbf{u}_j) \|_{[L^{s}(\Omega)]^4}^{r-1}\label{eqn:result_2},
\end{eqnarray}
where $C_i>0$ depends on $k_i$ whereas $C>0$ is independent of $i,j$. When $r>2$ it holds, instead,
 \begin{align}
    \nonumber &\| p_1 - p_2 \|_{L^{r'}(\Omega)} \\ &\leq  C_i
     \| k_1(| \varepsilon(\mathbf{u}_j) |) - k_2(| \varepsilon(\mathbf{u}_j) |)\|_{L^q(\Omega)}  
    \| \varepsilon(\mathbf{u}_j) \|_{[L^{s}(\Omega)]^4}\nonumber\\&\quad\times\left(1+\| \varepsilon(\mathbf{u}_1) \|_{[L^{r}(\Omega)]^4}^{\frac{r-2}{r}}+\| \varepsilon(\mathbf{u}_2) \|_{[L^{r}(\Omega)]^4}^\frac{r-2}{r}\right)\label{eqn:result_3},
\end{align}
where again $C_i>0$ depends on  $k_i$.
\end{theorem}
\begin{remark}
     In general it is not guaranteed that $\| k_1(| \varepsilon(\mathbf{u}_j) |) - k_2(| \varepsilon(\mathbf{u}_j) |)\|_{L^\infty(\Omega)}$ is finite, since there could be singularities in the tensor $\varepsilon(\mathbf{u}_j)$, leading to some explosion (in some particular case) in the quantities $k_i$, $i=1,2$. This means that we must find, if possible, a sufficiently small $q\geq r'$ such that $\| k_1(| \varepsilon(\mathbf{u}_j) |) - k_2(| \varepsilon(\mathbf{u}_j) |)\|_{L^q(\Omega)}<+\infty$. Clearly, if $q$ is too small, this could lead to the apparently unavoidable necessity of an extra assumption on the regularity of $\varepsilon(\mathbf{u}_j)$, which belongs \textit{a priori} only to $[L^r(\Omega)]^4$. For a more detailed theory ensuring more regularity on $\varepsilon(\mathbf{u}_j)$ we refer, for instance, to \cite{B1,B2, GPPV} and the references therein.
    \label{r>2}
\end{remark}
\begin{proof}
 Let us first prove \eqref{eqn:result}.   
 From the monotonicity of $A_i$ (cf. \cite[(2.19) and (2.14)]{barretnonnewton}), having chosen $A_i$ so that $r=\min_{l=1,2}\{r_l\}=r_i$, we have, for $j\not=i$, 
    \begin{eqnarray} 
       C_0 \| \mathbf{u}_1 - \mathbf{u}_2 \|_{[W^{1,r}(\Omega)]^2} &\leq & \left\langle A_i \mathbf{u}_i - A_i \mathbf{u}_j, \mathbf{u}_i - \mathbf{u}_j \right\rangle_{\mathbf{X}_i^*} \nonumber \\
                & =& \underbrace{\left\langle  A_i \mathbf{u}_i - A_j \mathbf{u}_j, \mathbf{u}_i - \mathbf{u}_j \right\rangle_{\mathbf{X}^*_i}}_{(i)}  +\underbrace{\left\langle  A_j \mathbf{u}_j - A_i \mathbf{u}_j, \mathbf{u}_i - \mathbf{u}_j \right\rangle_{\mathbf{X}^*_i}}_{(ii)}\nonumber.
        \end{eqnarray}
    First,  in view of $\nabla \cdot \mathbf{u}_i = 0$, we have
    \begin{equation*}
    \begin{split}
        (i) &= \left\langle        p_j, \nabla \cdot (\mathbf{u}_i - \mathbf{u}_j) \right\rangle  
        -\left\langle        p_i, \nabla \cdot (\mathbf{u}_i - \mathbf{u}_j) \right\rangle = 0.
    \end{split}
    \end{equation*} 
    On the other hand, for $q\in[r',\infty]$ such that $\| k_1(| \varepsilon(\mathbf{u}_j) |) - k_2(| \varepsilon(\mathbf{u}_j) |)\|_{L^q(\Omega)}<+\infty$,  we have
    \begin{eqnarray} 
            (ii) =\int_\Omega (k_j(| \varepsilon(\mathbf{u}_j) |)  - k_i(| \varepsilon(\mathbf{u}_j) |)    )\varepsilon(\mathbf{u}_j):\varepsilon(\mathbf{u}_i - \mathbf{u}_j), 
    \end{eqnarray} 
    which yields, after the use of the Hölder inequality, the following estimate
    \begin{eqnarray}  
            (ii)&\leq&\| k_1(| \varepsilon(\mathbf{u}_j) |) - k_2(| \varepsilon(\mathbf{u}_j) |)\|_{L^q(\Omega)} \| \varepsilon(\mathbf{u}_j) \|_{[L^{s}(\Omega)]^4} \| \varepsilon(\mathbf{u}_i - \mathbf{u}_j) \|_{[L^{r}(\Omega)]^4} 
        \nonumber          \nonumber \\&\leq& C \|k_1(|\varepsilon(\mathbf{u}_j)|)  - k_2(| \varepsilon(\mathbf{u}_j) |)\|_{L^q(\Omega)} \| \varepsilon(\mathbf{u}_j) \|_{[L^{s}(\Omega)]^4} \| \mathbf{u}_i - \mathbf{u}_j\|_{[W^{1,r}(\Omega)]^2},
    \end{eqnarray}
    where $s:=\frac{1}{1-\frac 1 q-\frac 1 r}$ if $q>r'$, otherwise, if $q=r'$, $s=\infty$. Here $C > 0$ is a generic constant. 
    %Notice that, when $r_i>2$ for $i=1,2$, we immediately deduce that $q=\infty$ and $s=r'$ (see Remark \ref{r>2}).  
    Combining the bounds for  $(i)$ and $(ii)$  we obtain:
    \begin{equation}
        \begin{aligned}
            \| \mathbf{u}_1 - \mathbf{u}_2 \|_{[W^{1,r}(\Omega)]^2} & \leq C \|k_1(|\varepsilon(\mathbf{u}_j)|)  - k_2(| \varepsilon(\mathbf{u}_j) |)\|_{L^q(\Omega)}\| \varepsilon(\mathbf{u}_j) \|_{[L^{s}(\Omega)]^4}.\\
        \end{aligned}
        \label{est_u}
    \end{equation}
 Concerning the pressure, let us first recall that, by the inf-sup condition \eqref{inf-sup:cont}, 
\begin{align}
\beta(r) \| p_1-p_2\|_{L^{r'}}&\leq \nonumber\sup_{\mathbf{w}\in W^{1,r}(\Omega)} \frac{\langle A_i \mathbf{u}_i - A_j \mathbf{u}_j, \mathbf{w}\rangle_{\mathbf{X}^*_i}}{\| \mathbf{w}\|_{[W^{1,r}(\Omega)]^2}}\\&\leq 
\sup_{\mathbf{w}\in [W^{1,r}(\Omega)]^2} \frac{\overbrace{\langle A_i \mathbf{u}_i - A_i \mathbf{u}_j, \mathbf{w}\rangle_{\mathbf{X}^*_i}}^{(iii)} - \overbrace{\langle A_j \mathbf{u}_j - A_i \mathbf{u}_j, \mathbf{w}\rangle_{\mathbf{X}^*_i}}^{(iv)}}{\| \mathbf{w}\|_{[W^{1,r}(\Omega)]^2}}\\&\leq 
\sup_{\mathbf{w}\in [W^{1,r}(\Omega)]^2} \frac{\vert (iii)\vert}{\| \mathbf{w}\|_{[W^{1,r}(\Omega)]^2}}
+\sup_{\mathbf{w}\in [W^{1,r}(\Omega)]^2} \frac{\vert (iv)\vert}{\| \mathbf{w}\|_{[W^{1,r}(\Omega)]^2}}
.
\end{align}
Now, the second summand is estimated exactly as the term $(ii)$ above, leading to
  \begin{eqnarray}  
            \vert(iv)\vert&\leq&\| k_1(| \varepsilon(\mathbf{u}_j) |) - k_2(| \varepsilon(\mathbf{u}_j) |)\|_{L^q(\Omega)} \| \varepsilon(\mathbf{u}_j) \|_{[L^{s}(\Omega)]^4} \| \varepsilon(\mathbf w) \|_{[L^{r}(\Omega)]^4} 
        \nonumber          \nonumber \\&\leq& C \|k_1(|\varepsilon(\mathbf{u}_j)|)  - k_2(| \varepsilon(\mathbf{u}_j) |)\|_{L^q(\Omega)} \| \varepsilon(\mathbf{u}_j) \|_{[L^{s}(\Omega)]^4} \| \mathbf{w}\|_{[W^{1,r}(\Omega)]^2},
    \end{eqnarray}
    entailing
  \begin{equation}
        \begin{aligned}
          \sup_{\mathbf{w}\in [W^{1,r}(\Omega)]^2} \frac{\vert (iv)\vert}{\| \mathbf{w}\|_{[W^{1,r}(\Omega)]^2}} & \leq C \|k_1(|\varepsilon(\mathbf{u}_j)|)  - k_2(| \varepsilon(\mathbf{u}_j) |)\|_{L^q(\Omega)}\| \varepsilon(\mathbf{u}_j) \|_{[L^{s}(\Omega)]^4},
        \end{aligned}
        \label{pp1}
    \end{equation}
    where $s,q$ are defined above.
Concerning $(iii)$, we need to distinguish between two cases.
First, we recall that by \cite[(2.1a), Lemma 2.1]{barretnonnewton}, it holds
\begin{align}
\label{basic}&\vert k_i(\vert \varepsilon(\mathbf{u}_i)\vert)\varepsilon(\mathbf{u}_i) -k_i(\vert \varepsilon(\mathbf{u}_j)\vert)\varepsilon(\mathbf{u}_j)\vert\nonumber\\&\leq C_i\vert \varepsilon(\mathbf u)\vert \left[(\vert\varepsilon(\mathbf{u}_1) \vert+\vert\varepsilon(\mathbf{u}_2)\vert)^\alpha(1+\vert\varepsilon(\mathbf{u}_1) \vert+\vert\varepsilon(\mathbf{u}_2)\vert)^{1-\alpha}\right]^{r-2},
\end{align}
where we set $\mathbf u:=\mathbf{u}_1-\mathbf{u}_2$.
Let now $r\leq 2$. Then, by recalling that $\frac 1 2(\vert x\vert +\vert y\vert)\leq \vert x\vert +\vert x-y\vert\leq 2(\vert x\vert +\vert y\vert)$, for any $x,y\in \mathbb{R}^{2\times 2}$, we obtain
\begin{align*}
 &\vert k_i(\vert \varepsilon(\mathbf{u}_i)\vert)\varepsilon(\mathbf{u}_i) -k_i(\vert \varepsilon(\mathbf{u}_j)\vert)\varepsilon(\mathbf{u}_j)\vert\\&\leq C_i\vert \varepsilon(\mathbf u)\vert \left[\frac{1}{2}(\vert\varepsilon(\mathbf{u}_1) \vert+\vert\varepsilon(\mathbf{u})\vert)^\alpha(2+\vert\varepsilon(\mathbf{u}_1) \vert+\vert\varepsilon(\mathbf{u})\vert)^{1-\alpha}\right]^{r-2}\\&\leq 
 C_i\vert \varepsilon(\mathbf u)\vert \frac{1}{2^{r-2}}(\vert\varepsilon(\mathbf{u}_1) \vert+\vert\varepsilon(\mathbf{u})\vert)^{r-2}\leq  \frac{C_i}{2^{r-2}}\vert \varepsilon(\mathbf u)\vert^{r-1}.
\end{align*}
Therefore, by H\"{o}lder's inequality, we deduce 
$$
\sup_{\mathbf{w}\in [W^{1,r}(\Omega)]^2} \frac{\vert (iii)\vert}{\| \mathbf{w}\|_{[W^{1,r}(\Omega)]^2}}\leq \frac{C_i}{2^{r-2}}\Vert \vert\varepsilon(\mathbf u)\vert^{r-1}\Vert_{L^{r'}(\Omega)}\leq \frac{C_i}{2^{r-2}}\Vert \varepsilon(\mathbf u)\Vert_{[L^{r}(\Omega)]^4}^{r-1}.
$$
Notice that $C_i$ depends on $k_i$.
Exploiting now \eqref{est_u} and putting together the estimate above and \eqref{pp1}, we immediately deduce \eqref{eqn:result_2}.
The case $r>2$ is much easier. The term involving $(iv)$ is estimated in \eqref{pp1}. Concerning the term involving $(iii)$, again from \eqref{basic} we infer, applying H\"{o}lder's inequality,
\begin{align*}
   & \sup_{\mathbf{w}\in W^{1,r}(\Omega)} \frac{\vert (iii)\vert}{\| \mathbf{w}\|_{[W^{1,r}(\Omega)]^2}}\\&\leq C_i\Vert \varepsilon(\mathbf u)\Vert_{[L^r(\Omega)]^4}\Vert (\vert\varepsilon(\mathbf{u}_1) \vert+\vert\varepsilon(\mathbf{u}_2)\vert)^\alpha(1+\vert\varepsilon(\mathbf{u}_1) \vert+\vert\varepsilon(\mathbf{u}_2)\vert)^{1-\alpha}\Vert_{L^r(\Omega)}^{\frac{r-2}{r}} \\&
\leq C_i\Vert \varepsilon(\mathbf u)\Vert_{[L^r(\Omega)]^4}\Vert 1+\vert\varepsilon(\mathbf{u}_1) \vert+\vert\varepsilon(\mathbf{u}_2)\vert\Vert_{[L^r(\Omega)]^4}^{\frac{r-2}{r}}\\&\leq
\tilde{C}_i\Vert \varepsilon(\mathbf u)\Vert_{[L^r(\Omega)]^4}(1+\Vert \varepsilon(\mathbf{u}_1) \Vert_{[L^r(\Omega)]^4}^{\frac{r-2}{r}}+\Vert\varepsilon(\mathbf{u}_2)\Vert_{[L^r(\Omega)]^4}^{\frac{r-2}{r}}),
\end{align*}
for some $\tilde{C}_i>0$, depending on $k_i$. Therefore, exploiting \eqref{basic} and putting together with \eqref{pp1} we easily end up with \eqref{eqn:result_3}. The proof is concluded.
\end{proof}

We now employ the above perturbation result to interpret the outcome of our campaign of numerical tests.  Let us first introduce ($\overline{\mathbf{u}}_{ICNN},\overline{p}_{ICNN}$), formally defined as: 
$$
\overline{\mathbf{u}}_{ICNN}:=\lim_{h \rightarrow 0} \mathbf{u}_{h,ICNN}, \quad \quad \quad
\overline{p}_{ICNN}:=\lim_{h \rightarrow 0} p_{h,ICNN} .
$$
In the sequel we employ Theorem \ref{prop:ourprop} with $k_1(t) = k(t)$ given by the Carreau law \eqref{eqn:carreau} with  parameters $k_\infty=0,k_0=2,\lambda=2$ and $n=r\in \{1.2, 1.6, 2, 2.4, 2.8\}$ and $k_2(t) = \text{ICNN}(t)$ representing its convex neural network approximation. In this context, it is reasonable to  assume that $k_1,k_2$ verify Assumptions ({\bf A}) for the same value of the parameter $r$. Moreover, we assume that $\| \varepsilon(\overline{\mathbf{u}}_{\text{ICNN}}) \|_{[L^{{r'}}(\Omega)]^4} < \infty $, and  $\|\varepsilon({\mathbf{\overline{u}}}) \|_{[L^{{r'}}(\Omega)]^4} < \infty$.

From Theorem \ref{prop:ourprop} with $q=\infty$, we immediately obtain
the following bounds: 
\begin{align}
    \nonumber\varepsilon_{\text{ICNN},u} &=\nonumber\| \overline{\mathbf{u}} - \overline{\mathbf{u}}_{ICNN} \| _{[W^{1,r}(\Omega)]^2} \\&\lesssim \| k(\vert \varepsilon(\overline{\mathbf{u}}_{ICNN})\vert) - \text{ICNN}(\vert \varepsilon(\overline{\mathbf{u}}_{ICNN})\vert)\|_{L^\infty(\Omega)}
    %\nonumber\\&
    \lesssim \| k - \text{ICNN} \|_{L^\infty(0,\infty)} \label{eq:aux:1},
    \end{align}
    \begin{align}
\varepsilon_{\text{ICNN},p} &= \| \overline{p} - \overline{{p}}_{\text{ICNN}} \| _{L^{r'}(\Omega)} \nonumber\\&\lesssim  \| k(\vert \varepsilon(\overline{\mathbf{u}}_{ICNN})\vert) - \text{ICNN}(\vert \varepsilon(\overline{\mathbf{u}}_{ICNN})\vert)\|_{L^\infty(\Omega)}^m
    %\\&
    \lesssim \| k - \text{ICNN}\|_{L^\infty(0,\infty)}^m\label{eq:aux:2}
  \end{align}
where $m=r-1$ if $r\in(1,2]$, $m=1$ if $r>2$, and the hidden constant may depend on the $W^{1,r'}$-norm of $\overline{\mathbf{u}}_{\ICNN}$ or  $\overline{{\bf u}}$, see Theorem \ref{prop:ourprop}. 
For the sake of conciseness we restrict ourselves to the case of the velocity and in  Table \ref{table:epsicnn} we report $\varepsilon_{\text{ICNN},u}$ (where the pair ($\overline{\mathbf{u}}_{ICNN},\overline{p}_{ICNN}$) is approximated by ($\mathbf{u}_{\overline{h},ICNN},p_{\overline{h},ICNN}$), with $\overline{h}=0.0027$) together with the error between the viscous law and the Input-Convex Neural Network measured in $L^\infty$-norm.

\begin{table}[H]
\centering
\begin{tabular}{|c| c| c |}
\hline
& $\varepsilon_{\text{ICNN},u}$ & $\| k - \text{ICNN}_{\bm{\theta}} \|_{L^\infty}$ \\
\hline 
$r=1.2$ & 0.00093 & 0.068426  \\
$r=1.6$ & 0.00051 & 0.004010  \\
$r=2.0$ & 0.00012 & 0.000550  \\
$r=2.4$ & 0.00067 & 0.005929 \\
$r=2.8$ & 0.00070 & 0.028900  \\
\hline
\end{tabular}
\caption{For different values of $r$, the behaviour of  $\varepsilon_{\text{ICNN},u}$ and   $L^\infty$-error are similar: they both decrease for increasing values of $r$ until they reach a minimum for $r=2$, then they both increase for $r>2$. 
}
\label{table:epsicnn}
\end{table}

Consider now the following simple inequality:
\begin{align}
        \nonumber
        \| \overline{\mathbf{u}} - \mathbf{u}_{h,ICNN}  \|_{[W^{1,r}(\Omega)]^2} 
         &\leq \|  \overline{\mathbf{u}}_{ICNN} - \mathbf{u}_{h,ICNN} \|_{[W^{1,r}(\Omega)]^2} + 
        \| \overline{\mathbf{u}} - \overline{\mathbf{u}}_{ICNN} \|_{[W^{1,r}(\Omega)]^2} \\ 
         & \lesssim h^{lr/2} + \| k - \text{ICNN} \|_{L^\infty(0,\infty)} 
         \label{eq:fem+ICNN}
\end{align}
where in the last step we employed \eqref{eqn:err_uh} and $l$ is the polynomial degree employed in the finite element approximation of the velocity $\mathbf{u}$.  A closer parallel inspection of  Figure \ref{fig:nn_errors} and Table \ref{table:epsicnn} reveals that the error $\| \overline{\mathbf{u}} - \mathbf{u}_{h,ICNN}  \|_{[W^{1,r}(\Omega)]^2}$ decreases with the expected rate $2$ (cf. Remark \ref{Rem:rate}),  until it reaches a plateau value corresponding  to $\| k - \text{ICNN} \|_{L^\infty(0,\infty)}$ (cf. the third column of Table \ref{table:epsicnn}). This shows that having a good approximation of $k$ through ICNN is essential to minimize the value of the plateau and thus the error.\\

We conclude this section addressing the convergence of ($\mathbf{u}_{h,ICNN},p_{h,ICNN}$) towards the solution  ($\overline{\mathbf{u}}_{ICNN},\overline{p}_{ICNN}$). Here, as above,  this latter is approximated with ($\mathbf{u}_{\overline{h},ICNN},p_{\overline{h},ICNN}$). The obtained convergence orders are reported in Table \ref{table:conv_up} and are in agreement with the quadratic theoretical optimal convergence \eqref{eqn:err_uh} (cf. Remark \ref{Rem:rate}).

\begin{table}[H]
\centering
\begin{tabular}{|c |c |c|}
\hline
& $\|\overline{\mathbf{u}}_{ICNN} - \mathbf{u}_{h,ICNN} \|_{[W^{1,r}(\Omega)]^2}$ & $\|\overline{p}_{ICNN} - p_{h,ICNN} \|_{L^{r'}(\Omega)}$ \\
\hline
$r=1.2$ & 2.02 & 1.80 \\
$r=1.6$ & 2.02 & 2.42 \\
$r=2.0$ & 2.02 & 3.14 \\
$r=2.4$ & 2.03 & 2.04 \\
$r=2.8$ & 2.04 & 2.37 \\
\hline
\end{tabular}
\caption{Convergence rates of $\mathbf{u}_{h,ICNN},p_{h,ICNN}$ towards the solution $\overline{\mathbf{u}}_{ICNN},\overline{p}_{ICNN}$. }
\label{table:conv_up}
\end{table}

\section{Conclusions}\label{S:Concl}
{
In this work, we introduced and tested on real-world scenarios a novel   theoretical framework and  computational strategy to build data-driven rheological models using Input-Convex Neural Networks. The results obtained demonstrate that the proposed approach represents a valuable and robust alternative to standard rheological laws. In the framework of the theoretical results of well-posedness for the non-Newtonian Stokes problem proved in \cite{barretnonnewton}, we have shown that the use of the data-driven ICNN rheological model is consistent with this theory. In particular, we demonstrated through numerical assessment that, by automatically choosing the viscosity law to be either convex or concave (depending on the flow regime), the mathematical conditions for ensuring the well-posedness of the differential problem are satisfied. 

Code and data that allow readers to reproduce the results in this paper are available at \url{https://github.com/Juu97/ICNN-Stokes}.
}

\section*{Funding}
MV and NP have been partially funded   by grant
PRIN2020 \emph{``Advanced polyhedral discretisations of heterogeneous PDEs for multiphysics problems''} n.20204LN5N5, funded by the Italian Ministry of Universities and Research (MUR).
The present research is part of the activities of ``Dipartimento di Eccellenza 2023-2027''. NP and MV are members of INdAM-GNCS. AP is member of INdAM-GNAMPA.

\bibliographystyle{siamplain}
\bibliography{icnn_arxiv}
\end{document}